\documentclass{elsart}
\usepackage{natbib}
\usepackage{amsmath,amssymb}
\usepackage[dvips]{graphicx}
\usepackage{psfrag}
\usepackage[all]{xy}

\newcommand{\I}{{\bf{I}}}
\newcommand{\Prime}{{\mathcal{P}}}
\newcommand{\Real}{{\Bbb R}}

\newcommand{\Zed}{{\Bbb Z}}
\newcommand{\Nat}{{\Bbb N}}

\newcommand{\upepsilon}{{ \xy
                   *{\xy (0,-3);(0,-2) **\dir{-} \endxy}
                   *\cir<3pt>{d^u}
                   \endxy}}

\newcommand{\Diff}{{\mathrm{Diff}}}
\newcommand{\PDiff}{{\mathrm{PDiff}}}
\newcommand{\FDiff}{{\mathrm{FDiff}}}
\newcommand{\PFDiff}{{\mathrm{PFDiff}}}
\newcommand{\KDiff}{{\mathrm{KDiff}}}
\newcommand{\dehn}{{\mathrm{DN}}}

\newcommand{\unknot}{{\mathcal{I}}}

\newcommand{\EK}[1]{{\mathrm{EC}({#1})}}
\newcommand{\EC}[2]{{{\mathrm{Emb}}(\Real^{#1} \times {#2},\Real^{#1} \times {#2})}}
\newcommand{\Emb}{{\mathrm{Emb}}}
\newcommand{\Imm}{{\mathrm{Imm}}}
\newcommand{\Cu}{{\mathcal{C}}}
\newcommand{\Aut}{{\mathrm{Aut}}}
\newcommand{\CAut}{{\mathrm{CAut}}}
\newcommand{\K}{{\mathcal{K}}}
\newcommand{\tK}{{\mathcal{\hat K}}}

\begin{document}
\runauthor{Budney}
\begin{frontmatter}
\title{Little cubes and long knots}
\author[Baddr]{Ryan Budney}

\address[Baddr]{IH\'ES, Le Bois-Marie, 35, route de Chartres, F-91440, Bures-sur-Yvette, France}
\address[Baddr]{Mathematics and Statistics, University of Victoria, PO BOX 3045 STN CSC, Victoria, B.C., Canada V8W 3P4}
\begin{abstract}
This paper gives a partial description of the
homotopy type of $\K$, the space of long knots in $\Real^3$.
The primary result is the construction of a homotopy equivalence
$\K \simeq \Cu_2(\Prime \sqcup \{*\})$
where $\Cu_2(\Prime \sqcup \{*\})$ is the free little $2$-cubes
object on the pointed space $\Prime \sqcup \{*\}$, where $\Prime \subset \K$
is the subspace of prime knots, and $*$ is a disjoint base-point.
In proving the freeness result, a close
correspondence is discovered between the Jaco-Shalen-Johannson decomposition of knot
complements and the little cubes action on $\K$. 
Beyond studying long knots in $\Real^3$ we show that for any compact
manifold $M$ the space of embeddings of $\Real^n \times M$  in
$\Real^n \times M$ with support in $\I^n \times M$  admits
an action of the operad of little $(n+1)$-cubes. 
If $M=D^k$ this embedding space is the space of
framed long $n$-knots in $\Real^{n+k}$, and the action of the little
cubes operad is an enrichment of the monoid structure given by the 
connected-sum operation.
\end{abstract}
\end{frontmatter}

\section{Introduction}\label{INTRODUCTION}

A theorem of Morlet's \cite{Mor} states  that the topological group
$\Diff(D^n)$ of boundary-fixing, smooth diffeomorphisms of the unit
 $n$-dimensional closed disc is homotopy equivalent to the $(n+1)$-fold loop space
$\Omega^{n+1}\left( PL_n/O_n \right)$.
Morlet's method did not involve the techniques invented by Boardman, 
Vogt and May \cite{BV, May1} for recognizing iterated loop spaces, little cubes actions.
This paper begins by defining little cubes operad actions
on spaces of diffeomorphisms and embeddings, thus making the loop space
structure explicit. In Theorem \ref{littlecthm} it's proved the embedding space
$$\EK{k,M} = \{ f \in \EC{k}{M}, supp(f) \subset \I^k \times M\}$$
admits an action of the operad of little $(k+1)$-cubes. Here the support of
$f$, $supp(f) = \overline{\{ x \in \Real^k \times M : f(x)\neq x \}}$
and $\I = [-1,1]$.

The case $k=1$ and $M=D^2$ is of primary interest in this paper
as $\EK{1,D^2}$ is the space of framed long knots in $\Real^3$.
In section \ref{longk} the structure of $\EK{1,D^2}$ as a little $2$-cubes
object is determined.  It is shown in Proposition \ref{modpro} that
the little $2$-cubes action on $\EK{1,D^2}$ restricts to a subspace $\tK$ which is
 homotopy equivalent to $\K$, the space of long knots in $\Real^3$. 
 Moreover it is shown that as little $2$-cubes objects,
$\EK{1,D^2} \simeq \tK \times \Zed$.
In Theorem \ref{freeness} it is shown that $\tK$ is a
free little $2$-cubes object on the subspace of prime long knots $\tK \simeq \Cu_2 \left( \Prime \sqcup \{*\} \right)$.
Theorems \ref{freeness} and \ref{littlecthm} are the main theorems of
this paper. 

\begin{figure}\label{fig1}
{
\psfrag{ek1d2}[tl][tl][1.2][0]{$f \in \EK{1,D^2}$}
\psfrag{p1}[tl][tl][0.8][0]{$1$}
\psfrag{m1}[tl][tl][0.8][0]{$-1$}
$$\includegraphics[width=10cm]{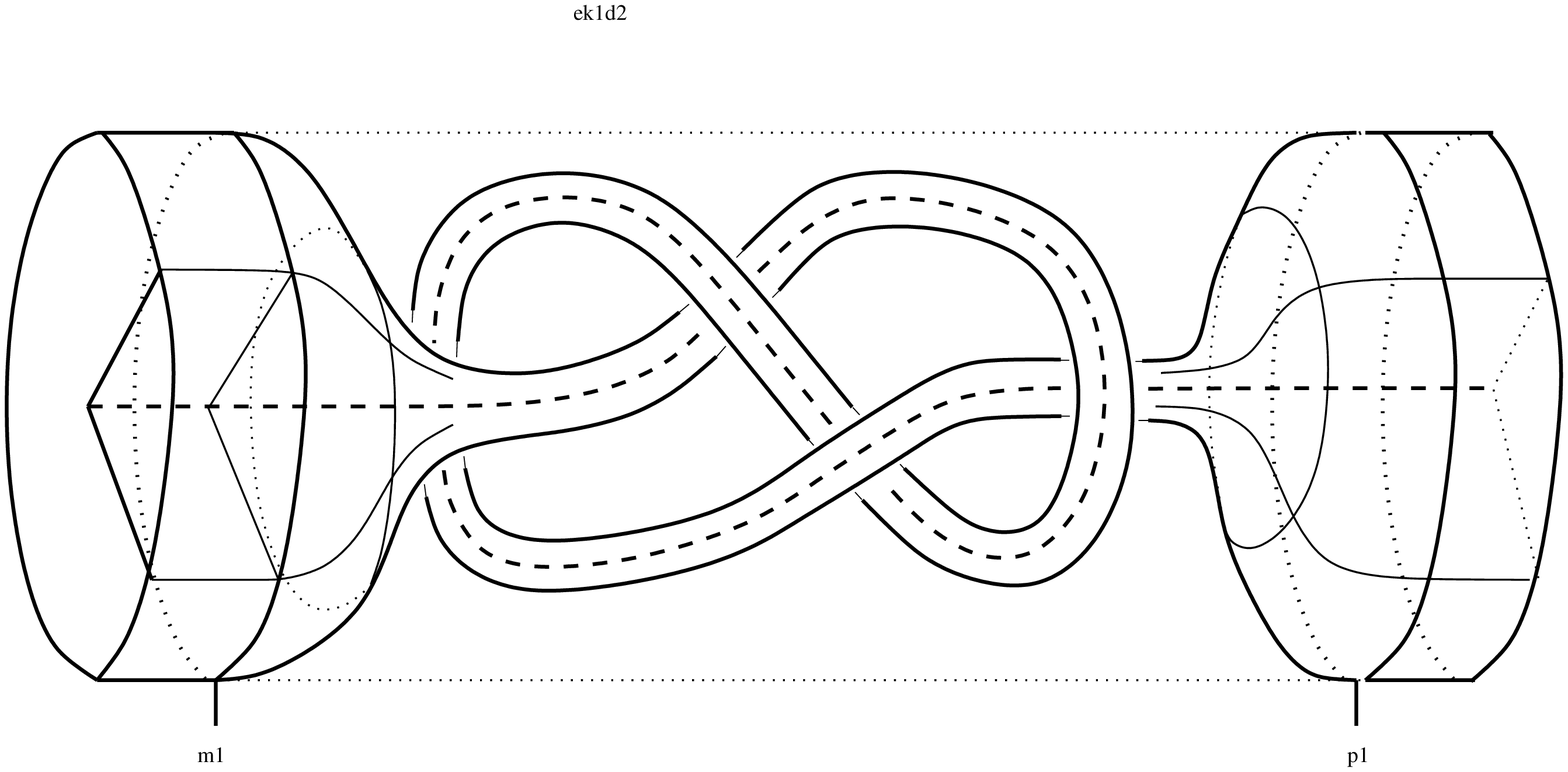}$$
\centerline{Figure 1}
}
\end{figure}

The homotopy-theoretic content of Theorem \ref{freeness} is that
$\K \simeq \Cu_2 \left( \Prime \sqcup \{*\} \right) \simeq \sqcup_{n=0}^{\infty} (\Cu_2(n) \times \Prime^n)/S_n$ where $\Cu_2(n)$ the space of $n$ little $2$-cubes. $\Cu_2(n)$ as an $S_n$-space
has the same homotopy type as the configuration space of $n$ labeled points in
the plane $C_n (\Real^2)$.  $\Prime \subset \K$
is the space of prime long knots, thus it is the union of all the components of $\K$ which
consist of prime knots. $S_n$ is the symmetric group on $n$ elements, acting
diagonally on the product. One interpretation of Theorem \ref{freeness} is that it
refines Schubert's Theorem \cite{Sch} which states that $\pi_0 \K$ is a free
commutative monoid with respect to the connected-sum operation
$\pi_0 \K \simeq \bigoplus_{\infty} \Nat$.  The refinement is a space-level
theorem about $\K$ where the cubes action on $\K$ replaces the connected-sum operation
on $\pi_0 \K$.
The novelty
of this interpretation is that the connected-sum is not a unique decomposition
in $\K$, as it is parametrized by a configuration space.  Perhaps the most
interesting aspect of Theorem \ref{freeness} is that it states
that the homotopy type of $\K$ is a functor
in the homotopy type of the space of prime long knots $\Prime$.
In Section \ref{endsec} we mention how the results in this paper
combine with results of Hatcher \cite{Hatcher4} and other results of
the author's \cite{topknot} to determine the full homotopy-type of $\K$.

There are elementary consequences of the little cubes actions
defined in Section \ref{littlec} that are of interest.
In Corollary \ref{cor1} we mention how the cubes action on $\EK{n,\{*\}}$
endows $\Diff(D^n) \simeq \EK{n,\{*\}}$ with the structure of
an $(n+1)$-fold loop space.  This corollary is part of Morlet's 
`Comparison' Theorem \cite{Mor}. To my knowledge, it is the 
first explicit demonstration of the $(n+1)$-cubes acting on
groups homotopy equivalent to $\Diff(D^n)$.
In Corollary \ref{cor2} the loop space recognition theorem together with
the cubes action on $\EK{k,D^m}$ and some elementary differential topology tell us
 that $\EK{k,D^m}$ is a $(k+1)$-fold loop space provided $m>2$. This
last result, to the best of my knowledge, is new.
Since these results appeared, Dev Sinha \cite{Dev2} has constructed an action
of the operad of $2$-cubes on the homotopy fiber of the map 
$\Emb(\Real,\Real^n) \to \Imm(\Real,\Real^n)$ for $n \geq 4$. Sinha's
result has recently been extended by Paolo Salvatore \cite{salvatore}, 
to construct actions of the operad of $2$-cubes on both the full embedding space $\Emb(\Real,\Real^n)$ and the `framed' long knot space $\EK{1,D^{n-1}}$ for 
$n \geq 4$, thus allowing for a comparison with the cubes actions 
constructed in this paper. Both the methods
of Salvatore and Sinha use the Goodwillie Calculus of Embeddings 
\cite{Good, Dev, DevKevin, Volic, Bud} together with the techniques of McClure and Smith \cite{McClure}. 

The existence of cubes actions on the space of long knots in $\Real^3$ was
conjectured by Turchin \cite{Tour}, who discovered a bracket on the
$E^2$-page of the Vassiliev spectral sequence  for the
homology of $\K$ \cite{Vass}.  Given the existence of a little $2$-cubes
action on $\EK{1,D^k}$ one might expect a co-bracket in the Chern-Simons approach
to the de Rham theory of spaces of knots \cite{Bott, Kont, Kohno, Cat1} but
at present only a co-multiplication is known \cite{Cat2}.  This paper could
also be viewed as an extension of the work of Gramain \cite{GramainPi1} who discovered
subgroups of the fundamental group of certain components of $\K$ which are
isomorphic to pure braid groups. 

\section{Actions of operads of little cubes on embedding spaces}\label{littlec}

In this section we define actions of operads of little cubes on various embedding spaces.
An invention of Peter May's, operads are designed to parametrize the
multiplicity of ways in which objects can be `multiplied'.  In the case of iterated loop spaces,
the relevant operad is the operad of little $n$-cubes, essentially defined by Boardman and Vogt \cite{BV} as `categories of operators in standard form,' and later recast into the 
language of operads by May \cite{May1}.

\begin{defn}
The space of long knots in $\Real^n$ is defined to be
$\Emb(\Real,\Real^n) = \{f : \Real \to \Real^n : \text{ where } f \text{ is a } C^\infty\text{-smooth embedding and } f(t)=(t,0,0,\cdots,0)$ for $|t|>1 \}$. 
We give $\Emb(\Real,\Real^n)$ the weak $C^\infty$ function space topology (see Hirsch \cite{Hirsch} \S 2.1).
$\Emb(\Real,\Real^n)$ is considered a pointed space with base-point 
given by $\unknot : \Real \to \Real^n$ where $\unknot(t)=(t,0,0,\cdots,0)$. 
Any knot isotopic to $\unknot$ is called an unknot. We reserve the
notation $\K$ for the space of long knots in $\Real^3$, ie: $\K = \Emb(\Real, \Real^3)$.
\end{defn}

The connected-sum operation $\#$ gives a homotopy-associative pairing
$$\# : \Emb(\Real,\Real^n) \times \Emb(\Real,\Real^n) \to \Emb(\Real,\Real^n)$$
As shown in Schubert's work \cite{Sch}, this pairing turns $\pi_0 \K$ (the path-components
of $\K$) into
a free commutative monoid with a countable number of generators (corresponding
to the isotopy classes of prime long knots). Schubert's
argument that $\pi_0 \K$ is commutative comes from the idea of
`pulling one knot through another,' illustrated in Figure 2.

\begin{figure}\label{fig2}
{
\psfrag{fg}[tl][tl][1][0]{$f \# g$}
\psfrag{gf}[tl][tl][1][0]{$g \# f$}
\psfrag{f}[tl][tl][1][0]{$f$}
\psfrag{g}[tl][tl][1][0]{$g$}
$$\includegraphics[width=12cm]{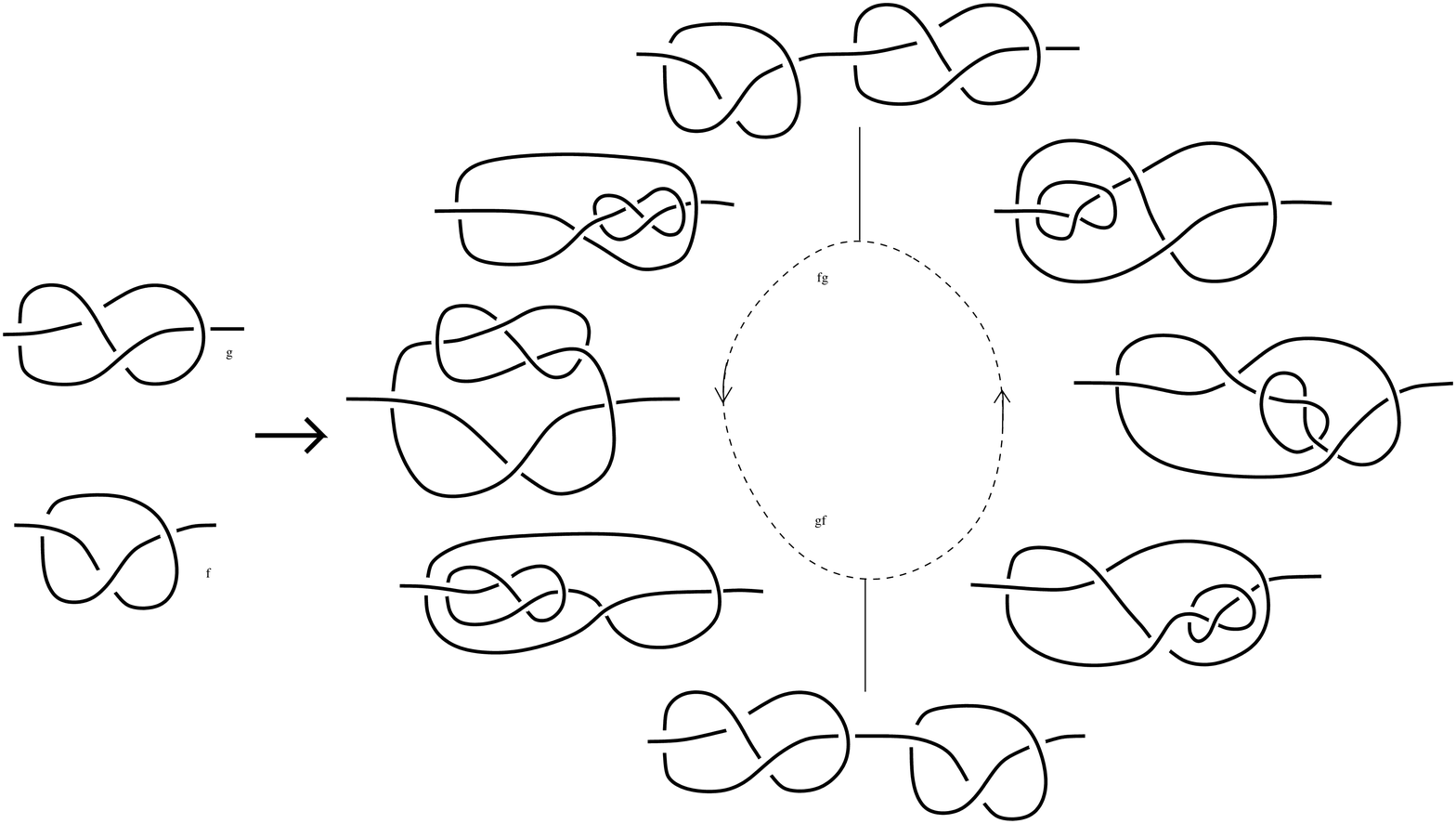}$$
\centerline{Figure 2}
}
\end{figure}

Figure 2 suggests the existence of a map $\iota: S^1 \times \K^2 \to \K$ such that $\iota(1,f,g) = f \# g$ and $\iota(-1,f,g)=g \# f$.  Such a map 
would exist if the connected sum operation on $\K$ was induced by a $2$-cubes action. 
Turchin's conjecture states that such a $2$-cubes action exists.

When first constructing the little $2$-cubes action on the
space of long knots, it was observed
that it is necessary to `fatten' the space $\K$ into a homotopy
equivalent space $\tK$ where the little cubes act.  The problem
with directly defining a little cubes action on $\K$ 
is that little cubes actions are very rigid.
Certain diagrams must commute \cite{May1, MSS}.  A homotopy commutative diagram
is not enough in the sense that one can not in general promote such diagrams to a
genuine cubes action. All known candidates for little cubes actions on $\K$
 that one might naively put forward have, at best, homotopy-commutative diagrams.
Definition \ref{firstbigdef} provides us a `knot space' $\EK{k,M}$ where
the connect-sum operation is given by composition of functions. The benefit
of this construction is that connect-sum becomes a strictly associative function,
allowing us to satisfy the rigid axioms of a cubes action.

\begin{figure}\label{fig3}
\centerline{$f \in \EK{1,D^2}$ and $L \in \CAut_1$} 
{
\psfrag{p1}[tl][tl][0.8][0]{$-1$}
\psfrag{m1}[tl][tl][0.8][0]{$1$}
\psfrag{mu}[tl][tl][1][0]{$\mu$}
\psfrag{e}[tl][tl][1][0]{$f$}
\psfrag{f}[tl][tl][1][0]{$L$}
\psfrag{e.f}[tl][tl][1][0]{$L.f$}
\psfrag{,}[tl][tl][1][0]{$,$}
$$\includegraphics[width=12cm]{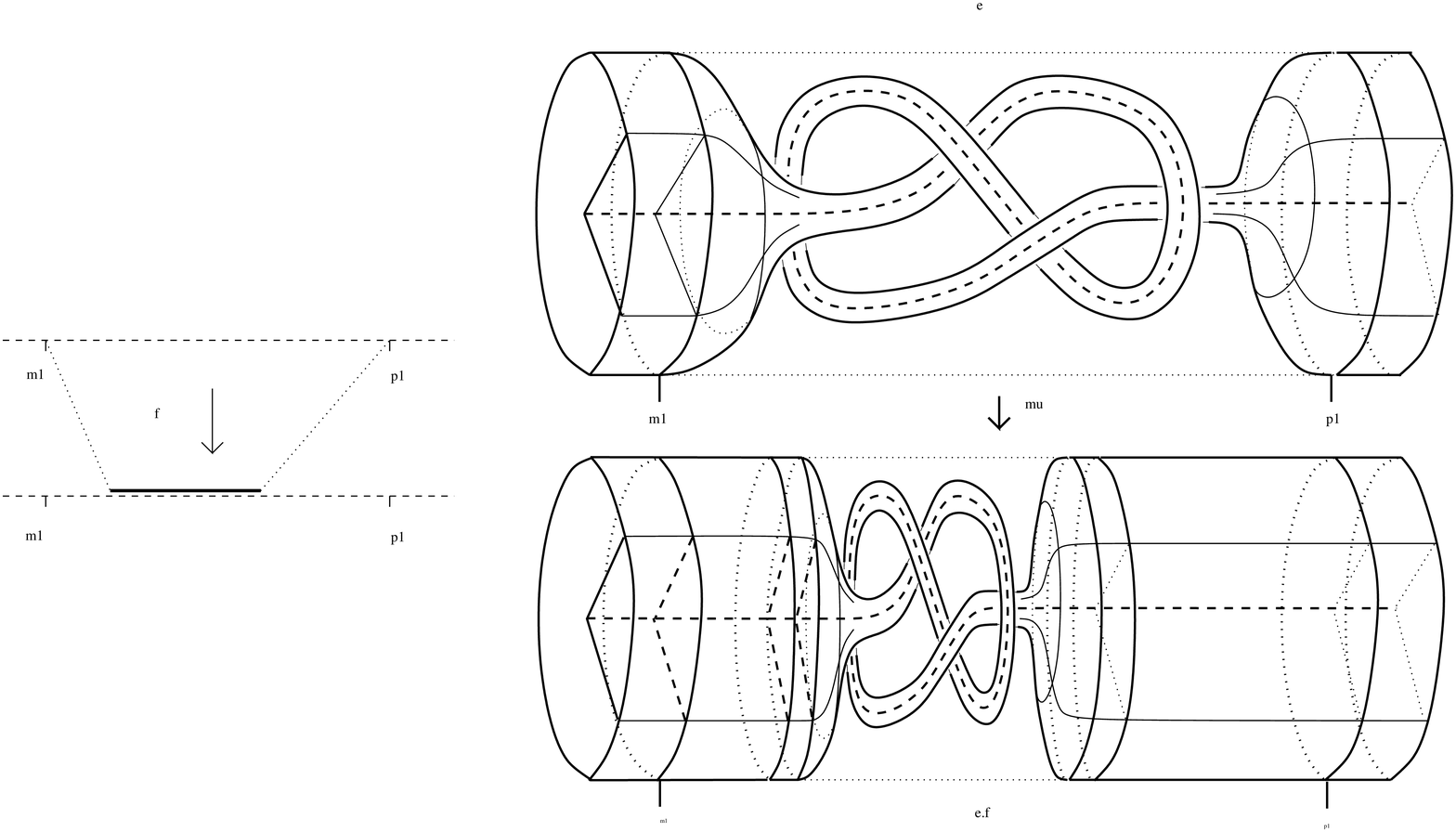}$$
\centerline{Figure 3}
}
\end{figure}

\begin{defn} \label{firstbigdef}
\begin{itemize}
\item $D^n := \{ x \in \Real^n : |x|\leq 1 \}$, where
 $\partial D^n = S^{n-1}$.
\item A (single) little $n$-cube is
a function $L : \I^n \to \I^n$ such that $L=l_1\times \cdots \times l_n$
where each $l_i : \I \to \I$ is affine-linear and increasing ie: $l_i(t)=a_i t + b_i$ for
some $a_i >0 $ and $b_i \in \Real$.
\item Let $\CAut_n$ denote the monoid of affine-linear
automorphisms of $\Real^n$ of the form $L=l_1 \times \cdots \times l_n$ where
$l_i$ is affine-linear and increasing for all $i \in \{1,2,\cdots,n\}$. 
\item Given a little
$n$-cube $L$, we sometimes abuse notation and consider $L \in \CAut_n$ by
taking the unique affine-linear extension of $L$ to $\Real^n$.
\item The space of $j$ little $k$-cubes  $\Cu_k(j)$ is the space of maps
$L : \sqcup_{i=1}^j \I^{k} \to \I^{k}$ such that
the restriction of
$L$ to the interior of its domain is an embedding,
and the restriction of
$L$ to any connected component of its domain
is a little $k$-cube.
Given $L \in \Cu_{k}(j)$, denote the restriction of $L$ to the
$i$-th copy of $\I^k$ by $L_i$. By convention $\Cu_k(0)$ is
taken to be a point. This makes the union $\sqcup_{j=0}^\infty \Cu_k(j)$
into an operad, called the operad of little $k$-cubes $\Cu_k$ \cite{May1, MSS}.
\item Given a compact manifold $M$, let $\EC{k}{M}$ denote the space of
$C^\infty$-smooth embeddings of $\Real^k \times M$ in $\Real^k \times M$. We do not
demand the embeddings to be proper ie: if $f \in \EC{k}{M}$ then the image of the
boundary of $\Real^k \times M$ need not lay in the boundary of $\Real^k \times M$.
We give this space the weak $C^\infty$-topology (See \cite{Hirsch} \S 2.1).
\item $\EK{k,M}$ is defined to be the subspace of $\EC{k}{M}$ consisting of
embeddings $f : \Real^k \times M \to \Real^k \times M$ whose support is contained
in $\I^k \times M$ ie: they are required to restrict to 
the identity function outside of $\I^k \times M$. We 
consider $\EK{k,M}$ to be a based space, with base-point given by the
identity function $Id_{\Real^k \times M}$. Any knot in the path component
of $Id_{\Real^k \times M}$ is typically called an unknot.
\end{itemize}
\end{defn}

We will show that the operad of little $(k+1)$-cubes acts on $\EK{k,M}$, but
first we define an action of the 
monoid $\CAut_k$ on $\EC{k}{M}$.
$$\mu : \CAut_k \times \EC{k}{M}  \to \EC{k}{M}$$
$$ \mu(L,f) = (L \times Id_M)\circ f \circ (L^{-1} \times Id_M) $$
In the above formula, we consider both $L$ and $L^{-1}$ to be
elements of $\CAut_n$. We
write the above action as $\mu(L,f)=L.f$ (see Figure 3).

\begin{prop}\label{contprop}
The two maps
$$\mu: \CAut_k \times \EC{k}{M} \to \EC{k}{M}$$
$$\circ : \EC{k}{M} \times \EC{k}{M} \to \EC{k}{M} $$
are continuous, where $\circ$ is composition.
\end{prop}

The continuity of  $\circ$ is an elementary consequence of the weak
topology.  The continuity of $\mu$ follows immediately.

\begin{defn}\label{littlecdef}
\begin{itemize}
\item Given $j$ little $(k+1)$-cubes, $L=(L_1,\cdots, L_j)\in \Cu_{k+1}(j)$
define the $j$-tuple of (non-disjoint) little $k$-cubes
$L^\pi = (L_1^\pi,\cdots, L_j^\pi)$ 
by the rule $L_i^\pi =l_{i,1} \times \cdots \times l_{i,k}$ where 
$L_i= l_{i,1} \times \cdots \times l_{i,k+1}$. Similarly define
$L^t \in \I^j$ by $L^t=(L_1^t,\cdots, L_j^t)$ where $L_i^t = l_{i,k+1}(-1)$
(see Figure 4).
\begin{figure}\label{fig4}
{
\psfrag{t}[tl][tl][1][0]{$\{0\}^n\times \Real$}
\psfrag{pit}[tl][tl][0.9][0]{$L^t$}
\psfrag{pic}[tl][tl][0.9][0]{$L^\pi$}
\psfrag{rn}[tl][tl][1][0]{$\Real^n\times\{0\}$}
\psfrag{f}[tl][tl][0.9][0]{$L$}
$$\includegraphics[width=5cm]{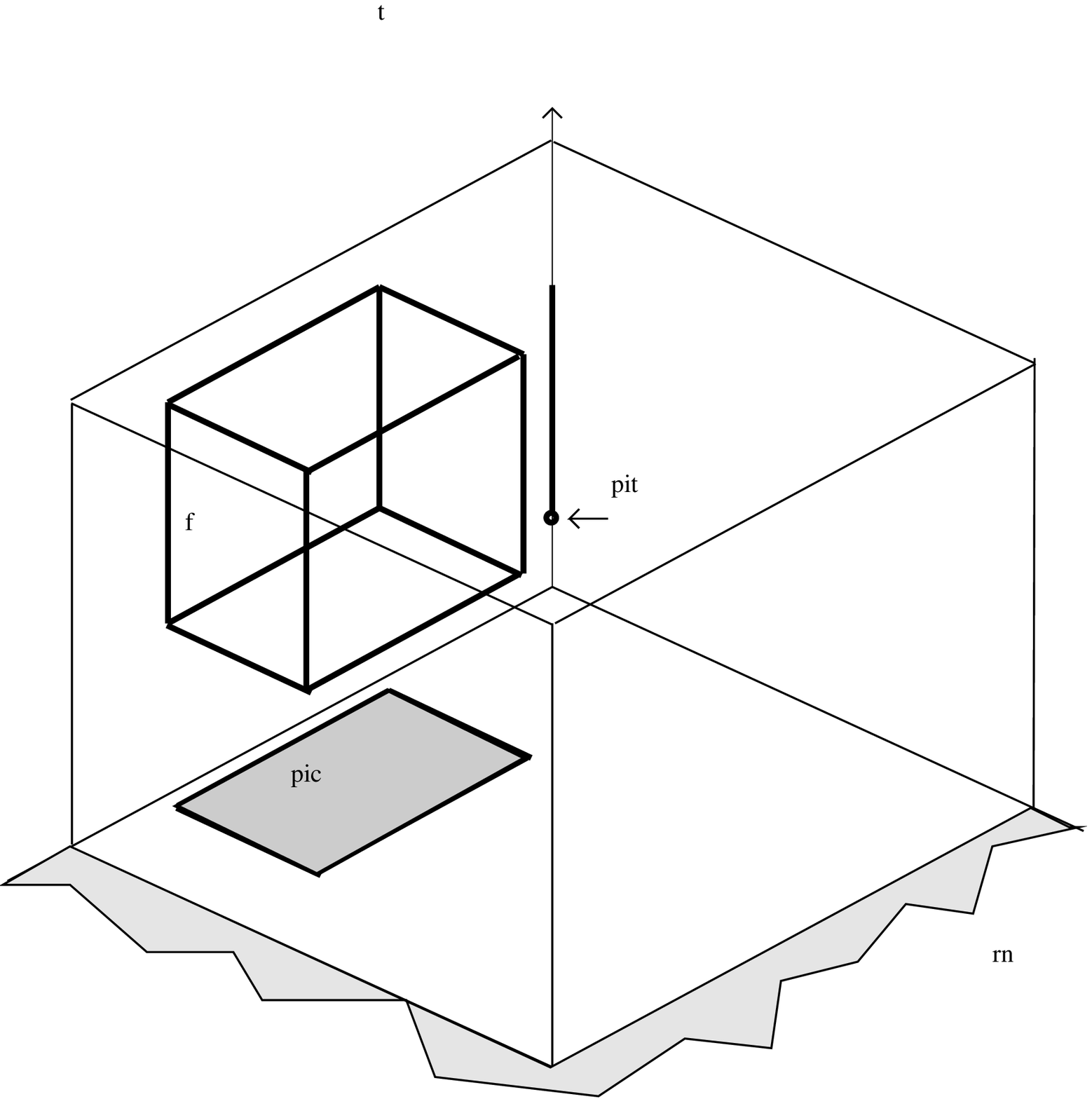}$$
\centerline{Figure 4}
}
\end{figure}
\item The action of the operad of little $(k+1)$-cubes on the space
$\EK{k,M}$ is given by
the maps $\kappa_j : \Cu_{k+1}(j) \times \EK{k,M}^j \to \EK{k,M}$
for $j \in \{1,2,\cdots\}$ defined by
$$\kappa_j(L_1,\cdots,L_j,f_1,\cdots,f_j) =
  L^\pi_{\sigma(1)}.f_{\sigma(1)}\circ L^\pi_{\sigma(2)}.f_{\sigma(2)}\circ\cdots \circ L^\pi_{\sigma(j)}.f_{\sigma(j)}$$
where $\sigma : \{1,\cdots, j\} \to \{1,\cdots,j\}$ is any permutation
such that $L^t_{\sigma(1)} \leq L^t_{\sigma(2)} \leq \cdots \leq L^t_{\sigma(j)}$.
The map $\kappa_0 : \Cu_{k+1}(0) \times \EK{k,M}^0 \to \EK{k,M}$
is the inclusion of a point $*$ in $\EK{k,M}$, defined
so that $\kappa_0(*) = Id_{\Real^k \times M}$ (see Figures 5 and 7).
\end{itemize}
\end{defn}

\begin{figure}\label{fig5}
\centerline{$L^t_1 < L^t_2$ so $\sigma$ is the identity and $\kappa_2(L_1,L_2,f_1,f_2)=L^\pi_1.f_1 \circ L^\pi_2.f_2$.}
{
\psfrag{1}[tl][tl][1][0]{$L_1$}
\psfrag{2}[tl][tl][1][0]{$L_2$}
\psfrag{p1}[tl][tl][0.7][0]{$1$}
\psfrag{m1}[tl][tl][0.7][0]{$-1$}
\psfrag{l1t}[tl][tl][1][0]{$L_1^t$}
\psfrag{l2t}[tl][tl][1][0]{$L_2^t$}
\psfrag{comma}[tl][tl][1][0]{,}
\psfrag{nu3}[tl][tl][1][0]{$\kappa_2$}
$$\includegraphics[width=13cm]{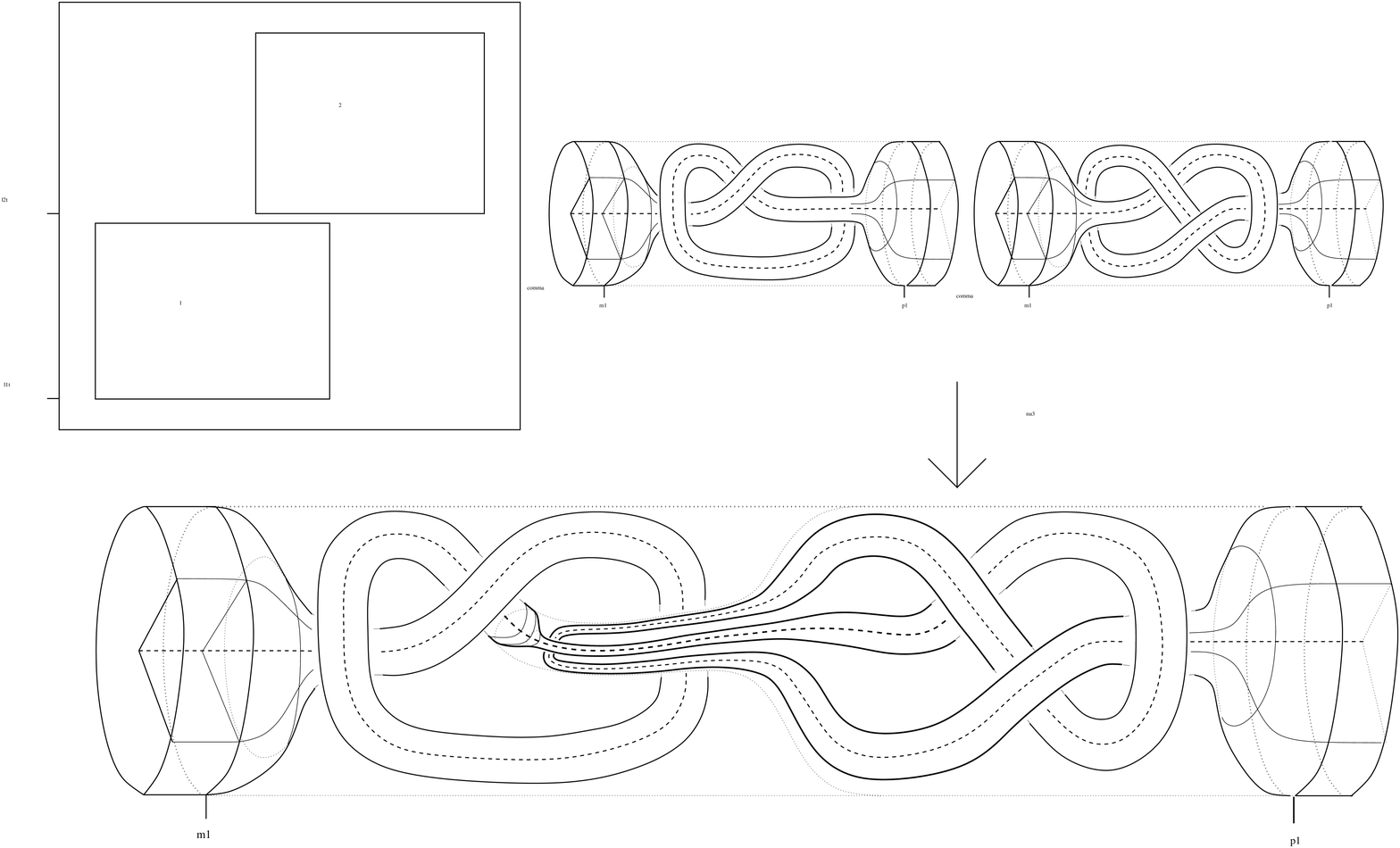}$$
\centerline{Figure 5}
}
\end{figure}

\begin{thm}\label{littlecthm}
For any compact manifold $M$ and any integer $k \geq 0$
the maps $\kappa_j$ for $j \in \{0, 1,2,\cdots\}$ define an action of the
operad of little $(k+1)$-cubes on $\EK{k,M}$.
\begin{pf}
First we show the map $\kappa_j$ is well-defined. The only ambiguity
in the definition is the choice of the permutation $\sigma$.
If there is an ambiguity in the choice of $\sigma$ this means that
a pair of coordinates $L^t_p$ and $L^t_q$ in $j$-tuple
$L^t=(L^t_1,\cdots,L^t_j)$ must be
equal. Since $L=(L_1,\cdots, L_j)$ are  disjoint cubes,
if a pair $L_p$ and $L_q$
have projections $L^t_p=L^t_q$, then $L^\pi_p$ and $L^\pi_q$ are disjoint.
Since $supp(L^\pi_p.f_p)=(L^\pi_p\times Id_M)(supp(f_p))$ and
$supp(L^\pi_{q}.f_q)=(L^\pi_{q}\times Id_M)(supp(f_q))$, $L^\pi_{p}.f_p$ and
$L^\pi_{q}.f_q$ must have disjoint support.
So the order of composition of $L^\pi_p.f_p$ and $L^\pi_q.f_q$ is irrelevant.
This proves the maps $\kappa_j$ are well-defined.

We prove the continuity of the maps $\kappa_j$.
Given a permutation $\sigma$ of the set
$\{1, \cdots, j\}$ consider the function
$$\kappa_\sigma : \Cu_{k+1}(j) \times \EK{k,M}^j \to \EK{k,M}$$ defined by
$$(L_1,\cdots,L_j,f_1,\cdots,f_j)
\longmapsto L^\pi_{\sigma(1)}.f_{\sigma(1)}\circ \cdots
 \circ L^\pi_{\sigma(j)}.f_{\sigma(j)}$$
This function is continuous, since the composition operation and the action
of  $\CAut_k$ is continuous by Proposition \ref{contprop}.   Given a
permutation $\sigma$, consider the subspace
$W_\sigma$ of $\Cu_{k+1}(j) \times \EK{k,M}^j$ where
$L^t_{\sigma(1)} \leq \cdots \leq L^t_{\sigma(j)}$.  Notice that our
map $\kappa_j$ when restricted to $W_\sigma$ agrees with
$\kappa_\sigma$. Thus the map $\kappa_j$ is the union of finitely
many continuous functions $\kappa_\sigma$ whose definitions agree where their
domains $W_\sigma$ overlap, so $\kappa_j$ is a continuous function
by the pasting lemma.

We need to show the maps $\kappa_j$ satisfy the axioms of a little
cubes action as described in sections 1 and 4 of \cite{May1}
(or II \S 1.4 of \cite{MSS}).
There are three conditions that must be satisfied: the
identity criterion, symmetry and associativity.  The identity
criterion is
tautological, since if $Id_{\I^{k+1}}$ is the identity little
$(k+1)$-cube, its projection is the identity cube, which acts
trivially on $\EK{k,M}$. Symmetry is similarly tautological.
The associativity condition demands that the diagram in Figure 6
commutes.
\begin{figure}\label{fig6}
\psfrag{bp2}[tl][tl][0.9][0]{$\Cu_{k+1}(n) \hskip 2mm \times$}
\psfrag{tr}[tl][tl][0.9][0]{$\Cu_{k+1}(n) \times \EK{k,M}^n$}
\psfrag{br}[tl][tl][0.9][0]{$\EK{k,M}$}
\psfrag{bl}[tl][tl][0.9][0]{$\Cu_{k+1}(j_1+\cdots+j_n) \times \EK{k,M}^{j_1+\cdots+j_n}$}
\psfrag{bp1}[tl][tl][0.9][0]{$\Cu_{k+1}(j_1) \times \EK{k,M}^{j_1} \times \cdots \times \Cu_{k+1}(j_n) \times \EK{k,M}^{j_n}$}
$$\includegraphics[width=12cm]{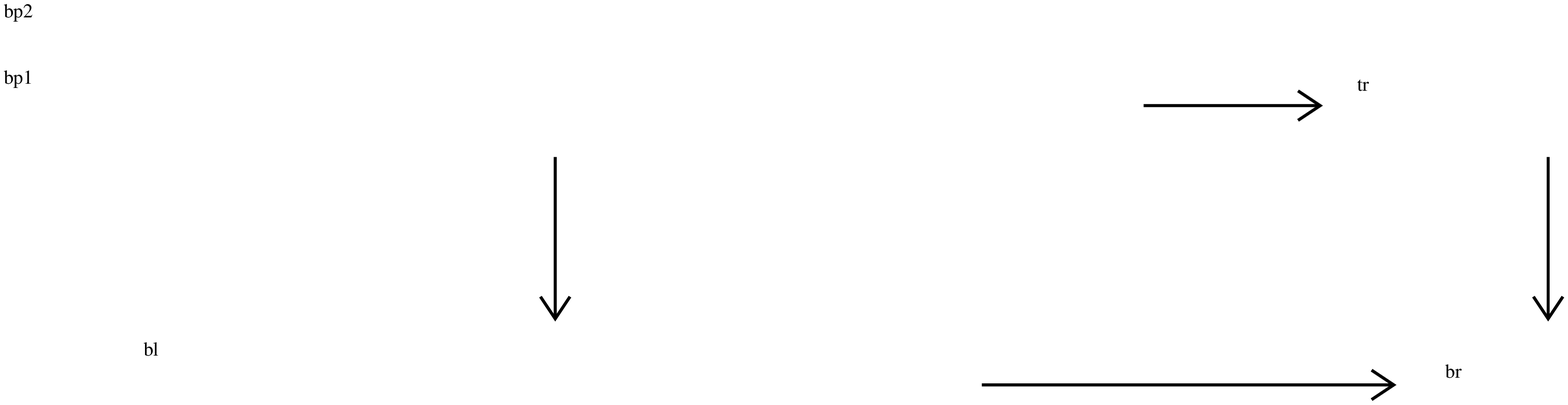}$$
\centerline{Figure 6}
\end{figure}
The commutativity of this diagram follows from the same argument
given that shows that the
maps are well-defined.  If one chases the arrows around the diagram
both ways,
the two objects that you get in $\EK{k,M}$ are composites of the same
embeddings, perhaps in a different order. 
Any pair of embeddings
that have their order permuted must have disjoint supports,
so the change in order of composition is irrelevant.
\end{pf}
\end{thm}

\begin{figure}\label{fig7}
$L^t_1 < L^t_3 < L^t_2$ so $\sigma=(23)$ and $\kappa_3(L_1,L_2,L_3,f_1,f_2,f_3)=L^\pi_1.f_1\circ L^\pi_3.f_3 \circ L^\pi_2.f_2$, which explains why we see the figure-8 knot `inside' of the trefoil.
{
\psfrag{1}[tl][tl][1][0]{$L_1$}
\psfrag{2}[tl][tl][0.7][0]{$L_2$}
\psfrag{3}[tl][tl][1][0]{$L_3$}
\psfrag{p1}[tl][tl][0.7][0]{$1$}
\psfrag{m1}[tl][tl][0.7][0]{$-1$}
\psfrag{comma}[tl][tl][1][0]{,}
\psfrag{nu3}[tl][tl][1][0]{$\kappa_3$}
\psfrag{l3t}[tl][tl][1][0]{$L_3^t$}
\psfrag{l2t}[tl][tl][1][0]{$L_2^t$}
\psfrag{l1t}[tl][tl][1][0]{$L_1^t$}
\psfrag{f3}[tl][tl][1][0]{$f_3$}
\psfrag{phi1}[tl][tl][1][0]{$f_1$}
\psfrag{phi2}[tl][tl][1][0]{$f_2$}
\psfrag{phi3}[tl][tl][1][0]{$f_3$}
$$\includegraphics[width=13cm]{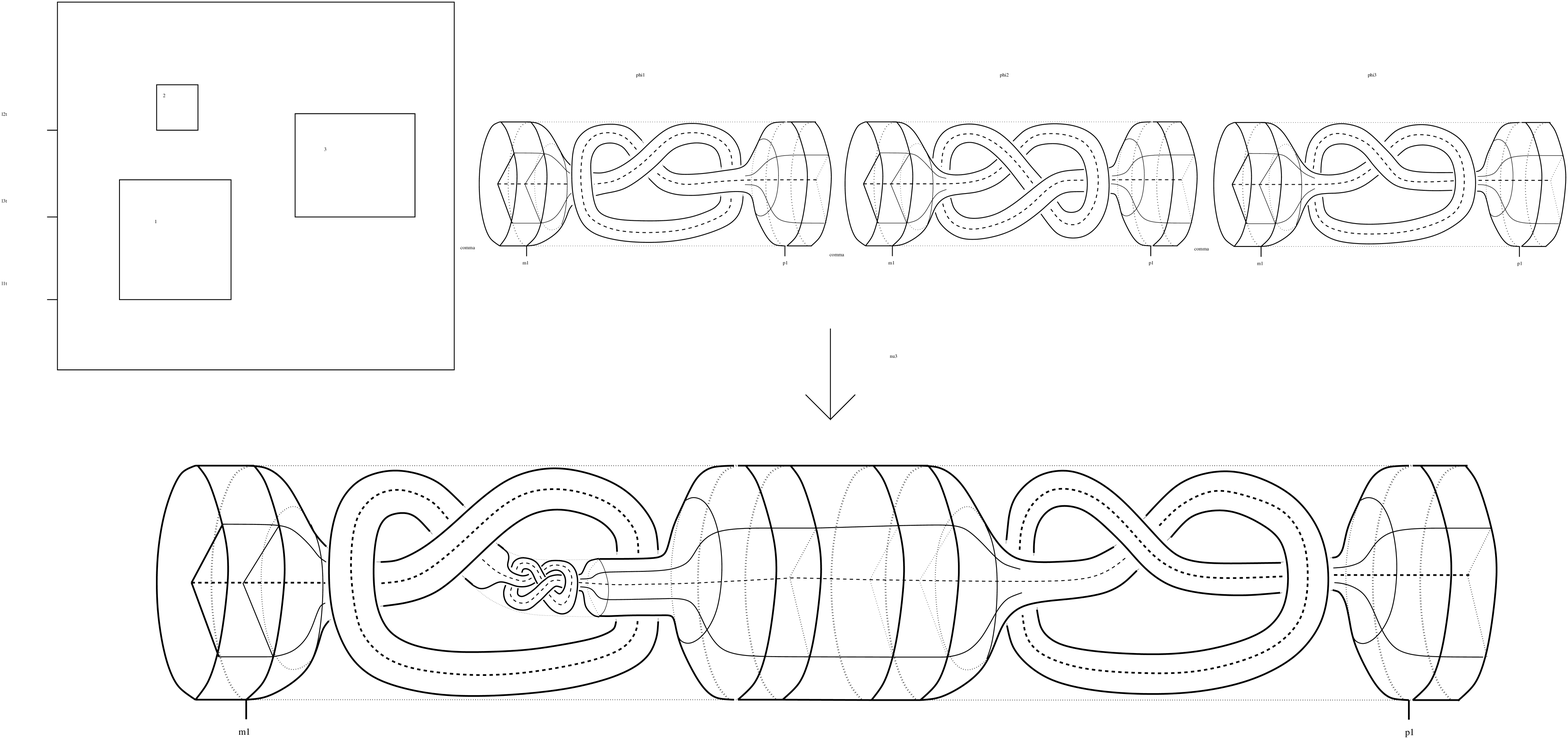}$$
\centerline{Figure 7}
}
\end{figure}

\begin{cor}\label{cor1}
The group of boundary-fixing diffeomorphisms of the compact
$n$-dimensional ball,
$\Diff(D^n)$ is homotopy-equivalent to an $(n+1)$-fold loop space.
\begin{pf}
Peter May's loop space recognition theorem \cite{May2} states that
 a little $(n+1)$-cubes
object $X$ is (weakly) homotopy equivalent to an $(n+1)$-fold loop space
if and only if the induced monoid structure on $\pi_0 X$ is a group.

Consider the monoid structure on $\pi_0 \EK{n,\{*\}}$.
Let $L=(L_1,L_2) \in \Cu_{n+1}(2)$ be two little $(n+1)$-cubes
such that $L^\pi=(L^\pi_1,L^\pi_2)=(Id_{\I^n},Id_{\I^n})$.
Suppose $L^t=(L^t_1,L^t_2)$ with $L^t_1 < L^t_2$, then
$\kappa_2(L_1,L_2,f_1,f_2)=f_1 \circ f_2$.  This
means the induced monoid structure on $\pi_0 \EK{n,\{*\}}$ is given by
composition. $\EK{n,\{*\}}$ is a group under composition
since it is the group of diffeomorphisms $\Real^n$ with support
contained in $\I^n$. Thus, $\pi_0 \EK{n,\{*\}}$ is also a group, and
so $\EK{n,\{*\}}$ is weakly homotopy equivalent to
an $(n+1)$-fold loop space. Since $\EK{n,\{*\}}$ has
the weak (compact-open) topology (see Hirsch \cite{Hirsch} \S 2.1)  it
satisfies the first axiom of countability and so
the topology on $\EK{n,\{*\}}$ is compactly-generated in the sense of Steenrod
\cite{Steenrod}.  Thus by the loop space recognition theorem,
the $\EK{n,\{*\}}$ is homotopy-equivalent to an $(n+1)$-fold
loop space.

 Provided we show that
$\Diff(D^n) \simeq \EK{n,\{*\}}$ we are done.
Fix a collar neighborhood of $S^{n-1}$. There is a restriction map from
$\Diff(D^n)$ to the space of collar neighborhoods of $S^{n-1}$ in $D^n$.
This restriction map is a fibration \cite{Pal} and the space of
collar neighborhoods of $S^{n-1}$ in $D^n$ is contractible (see \cite{Hirsch} \S 4.5.3).
The above argument is not sufficient, because the fiber of this
fibration is not
$\EK{n,\{*\}}$. Replace the smooth collar neighborhood of $S^{n-1}$ in $D^n$
with a manifold-with-corners neighborhood of $S^{n-1}$ 
which is the complement of an open cube in $D^n$.
In this case we get a fibration whose fiber we can identify with $\EK{n,*}$.
The argument that the space of cubical collar neighborhoods is contractible is analogous
to the proof in Hirsch's text (see \cite{Hirsch} \S 4.5.3).
\end{pf}
\end{cor}

May's recognition theorem applies equally-well to spaces that have actions
of the operad of (unframed) little balls \cite{MSS}.  Thus we could have simply adapted
Definition \ref{littlecdef} to give an action of the space of unframed
$(n+1)$-balls directly on the the space $\Diff(D^n)$ and deduced the
result without recourse to the intermediate homotopy-equivalence 
$\Diff(D^n) \simeq \Diff(\I^n)$.

The above corollary is also a corollary of Morlet's `Comparison Theorem' \cite{Mor}.
Morlet's manuscript was not widely distributed. A proof of Morlet's Theorem can
be found in Burghelea and Lashof's paper \cite{BL}, as well as in Kirby and
Siebenmann's book \cite{KS}. As Siebenmann points out, the Morlet Comparison
Theorem was first observed by Cerf.

\begin{cor}\label{cor2}
$\EK{k,D^n}$ is homotopy equivalent to a $(k+1)$-fold loop space, provided
$n>2$.
\begin{pf}
This follows from the loop space recognition theorem \cite{May2}
since we will show that $\pi_0 \EK{k,D^n}$ is a group. Consider 
the fibration $\EK{k,D^n} \to \Emb(\Real^k,\Real^{k+n})$
where $\Emb(\Real^k,\Real^{k+n})$ is the space
$\{ f : \Real^k \to \Real^{k+n}$ \-:\- $f(t_1,t_2,$ $\cdots,t_k)=(t_1,t_2,\cdots,t_k,0,0,\cdots,0)$  if 
$|t_i| \geq 1$ for any $i \in \{1,2,\cdots,k\} \}$.
Haefliger proved \cite{Haefliger2} that $\pi_0 \Emb(\Real^k,\Real^{k+n})$ 
is a group provided $n>2$ where the group structure is induced by
concatenation, thus $\pi_0 \EK{k,D^n}$ is a group as it is a monoid
which is an extension of two groups (see \cite{family} for an alternative
proof of Haefliger's theorem).
\end{pf}
\end{cor}

Our preferred model for $\K$
will be a subspace $\tK$ of $\EK{1,D^2}$,
which we will relate back to the standard model $\K$. Given an embedding $f \in \EK{1,D^2}$,
define $\omega(f) \in \Zed$ to be the linking number of 
$f_{|\Real \times \{(0,0)\}}$ with
$f_{|\Real \times \{(0,1)\}}$. One concrete way to define this integer is
as the transverse intersection number of the map 
$\Real^2 \ni (t_1,t_2) \longmapsto f(t_1,0,1)-f(t_2,0,0) \in 
 \Real^3 - \{(0,0,0)\}$ with
 the ray $\{(0,t,0) : t>0\} \subset \Real^3 - \{(0,0,0)\}$.
$\omega(f)$ is called the framing number of $f$.
$\omega : \EK{1,D^2} \to \Zed$ is a $2$-cubes equivariant fibration, 
and the
framing number is additive $\omega(f_1\circ f_2)=\omega(f_1)+\omega(f_2)$.
We consider $\Zed$ to be an abelian group, and thus a little $2$-cubes object.

\begin{defn}\label{fatdef} $\tK$, the space of `fat' long knots in $\Real^3$ is defined to be the kernel of $\omega$, $\tK = \omega^{-1}\{0\}$.
\end{defn}

\begin{prop}\label{modpro} The two spaces $\tK$ and $\K$ are homotopy equivalent.
\begin{pf}
Consider the fibration $\EK{1,D^2} \to \Emb(\Real,\Real^3)$ given by
restriction $f \longmapsto f_{|\Real \times \{(0,0)\}}$ \cite{Pal}.
Let $X$ denote the fiber of this fibration. By definition, $X$ is the space of tubular
neighborhoods of the unknot which are standard outside of $\I \times D^2$.
By the classification of tubular neighborhoods theorem
(see for example \cite{Hirsch} \S 4.5.3), $X$ is homotopy equivalent to
the space of fibrewise-linear automorphisms of $\Real \times D^2$ with
support in $\I \times D^2$, ie: $X \simeq \Omega SO_2 \simeq \Zed$.
Thus $\omega$ defines a splitting of the
fibration $X \to \EK{1,D^2} \to \Emb(\Real,\Real^3)$,
giving the two homotopy equivalences
$$\EK{1,D^2} \simeq \Emb(\Real,\Real^3) \times \Zed \hskip 10mm \xymatrix@R=5pt{\tK \ar@{}[d]|{\upepsilon} \ar[r]^-{\simeq} &
\K \ar@{}[d]|{\upepsilon} \\
 f \ar@{|->}[r] & f_{|\Real \times \{(0,0)\}} }$$
\end{pf}
\end{prop} 

Combining Proposition \ref{modpro} with the proof of Corollary \ref{cor2}
we get the following observation.

\begin{cor}There is an action of the operad of $(k+1)$-cubes on spaces
homotopy-equivalent to the `long embedding spaces'
$\Emb(\Real^k,\Real^{k+n})$ for all $k \in \Nat$ and $n \leq 2$.
\end{cor}

As mentioned in the introduction, Salvatore \cite{salvatore} has removed
the bound $n \leq 2$ in the above corollary, provided $k=1$.

\section{The freeness of the $2$-cubes action on $\tK$}\label{longk}

The goal of this section is to prove that $\tK \simeq \Cu_2(\Prime \sqcup \{*\})$,
where $\Prime \subset  \tK$ is the subspace of
prime knots. $\Prime = \{ f \in \tK : f$ 
is nontrivial and not a connected-sum of 2 or more
nontrivial knots$\}$.

If $X$ is a pointed space with base-point
$* \in X$ the free little $2$-cubes object on $X$ \cite{May1} is the space
$\Cu_2(X) = \left( \left(\sqcup_{n=0}^\infty \Cu_2(n) \times X^n\right)\kern-0.3em/S_n\right)\kern-0.3em/\kern-0.5em\sim$.
$S_n$ is the symmetric group, acting diagonally on the product 
in the standard way, and the equivalence
relation $\sim$ is generated by the relations
$$\left((f_1,\cdots,f_{i-1},f_i,f_{i+1},\cdots,f_n),(x_1,\cdots,x_{i-1},*,x_{i+1},\cdots,x_n)\right)$$
$$ \sim ((f_1,\cdots,f_{i-1},f_{i+1},\cdots,f_n),(x_1,\cdots,x_{i-1},x_{i+1},\cdots,x_n))$$
If we give an arbitrary unpointed space $X$ a disjoint base-point $*$, then
there is the identity 
$\Cu_2(X \sqcup \{*\})\equiv \sqcup_{n=0}^\infty (\Cu_2(n) \times X^n)/S_n$.
Thus, we will prove 
$\tK \simeq \sqcup_{n=0}^\infty \Cu_2(n) \times_{S_n} \Prime^n$.

\begin{thm}\label{freeness} $\tK \simeq \Cu_2(\Prime \sqcup \{*\})$, moreover 
the map $\sqcup_{n=0}^\infty \kappa_n : \sqcup_{n=0}^\infty \Cu_2(n) \times_{S_n} \tK^n \to \tK$ restricts to a homotopy equivalence
$$ \sqcup_{n=0}^\infty \Cu_2(n) \times_{S_n} \Prime^n \to \tK$$
\end{thm}

To prove Theorem \ref{freeness} we first build up a close correspondence between the
little cubes action and the satellite decomposition of
knots, or to be more precise, the JSJ-decomposition \cite{JacoShalen} 
of knot complements (also sometimes also known as the splice decomposition \cite{EN}) . We then 
use techniques of Hatcher's to reduce the proof of Theorem \ref{freeness} to a
problem about a diagram of mapping class groups of $2$ and $3$-dimensional manifolds.  
 
\begin{defn}\label{sumdef}
\begin{itemize}
\item Given a long knot $f \in \tK$, we denote the component of
$\tK$ containing $f$ by $\tK_f$. 
\item We say $f$ is a connected-sum of $f_1, \cdots, f_n$ if there
exists $\tilde f \in \tK_f$ with
 $\tilde f = \kappa_n (L_1,L_2,\cdots,L_n,f_1,f_2,\cdots,f_n)$,
for some $n$, $(L_1,L_2,\cdots,L_n)\in \Cu_2(n)$ and
$f_i \in \tK$ for all $i \in \{1,2,\cdots,n\}$. 
Denote this by $f \sim f_1 \# f_2 \# \cdots \# f_n$ and call the
long knots $\{f_i : i \in \{1,2,\cdots,n\}\}$ summands of $f$.
\item For any long knot $f \sim f \# Id_{\Real\times D^2}$.
If $f \sim f_1 \# f_2 \# \cdots \# f_n$ we call the connected-sum
trivial if $(n-1)$ of the long knots $\{f_1,f_2,\cdots,f_n\}$
are in $\tK_{Id_{\Real \times D^2}}$.
A long knot is prime if is not in the component of the unknot, and if all
connected-sum decompositions of it are trivial. 
\end{itemize}

Let $Q_i$ denote the $2$-cube 
$[-1+\frac{4i-2}{2n+1},-1+\frac{4i}{2n+1}]\times [0,\frac{2}{2n+1}]$.
Choose the base-point $*$ for $\Cu_2(n)$,
$*=(Q_1,Q_2,\cdots,Q_n)$ as in Figure 8.
\begin{figure}\label{fig8}
{
\psfrag{-1}[tl][tl][1][0]{$-1$}
\psfrag{i1}[tl][tl][0.7][0]{$-1+\frac{2}{2n+1}$}
\psfrag{i2}[tl][tl][0.7][0]{$-1+\frac{4}{2n+1}$}
\psfrag{i3}[tl][tl][0.7][0]{$-1+\frac{6}{2n+1}$}
\psfrag{i4}[tl][tl][0.7][0]{$-1+\frac{8}{2n+1}$}
\psfrag{im2}[tl][tl][0.7][0]{$-1+\frac{4n-2}{2n+1}$}
\psfrag{im1}[tl][tl][0.7][0]{$-1+\frac{4n}{2n+1}$}
\psfrag{1}[tl][tl][1][0]{$1$}
\psfrag{q1}[tl][tl][1][0]{$Q_1$}
\psfrag{q2}[tl][tl][1][0]{$Q_2$}
\psfrag{qn}[tl][tl][1][0]{$Q_n$}
\psfrag{yc}[tl][tl][1][0]{$\frac{2}{2n+1}$}
\psfrag{ddd}[tl][tl][1][0]{$\cdots$}
\psfrag{rx0}[tl][tl][1][0]{$\Real \times \{0\}$}
$$\includegraphics[width=13cm]{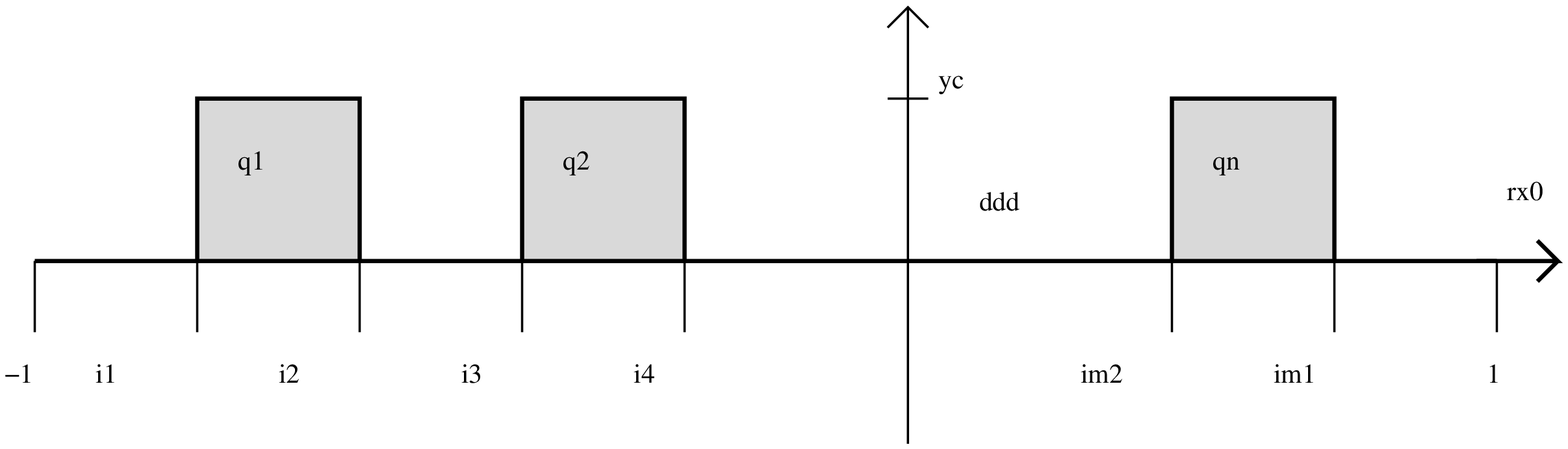}$$
\centerline{Figure 8}
}
\end{figure}

Since $\Cu_2(n)$ is connected, we can choose the $n$ little $2$-cubes
 $(L_1,L_2,\cdots,L_n)\in \Cu_2(n)$ in Definition \ref{sumdef} to be
 $*$.
As in Proposition \ref{modpro} we can associate to $f \in \tK$
the long knot $g \in \Emb(\Real,\Real^3)$ where $g=f_{|\Real \times \{(0,0)\}}$. Define
$B_i = \{ x \in \Real^3 : |x-(\frac{4i-2n-2}{2n+1},0,0)| \leq \frac{1}{2n+1} \}$
and $S_i = \partial B_i$ (see Figure 9).
 We will provide an equivalent definition
for $f$ to be a connected-sum in terms of $g$ (see Figure 9).

We say $g$ is a connected-sum if $g$ is
isotopic to $g' \in \Emb(\Real,\Real^3)$ such
that:

\begin{itemize}
\item $supp(g')\subset (\sqcup_{i=1}^n B_i)\cap (\Real\times \{0\}^2)$
\item $img(g') \cap S_i = (\Real \times \{0\}^2 )\cap S_i$ for all $i \in \{1,2,\cdots,n\}$.
\item There exists long knots (the summands of $g$)
$g_i \in \Emb(\Real,\Real^3)$ for $i \in \{1,2,\cdots,n\}$ such
that $supp(g_i) \subset B_i \cap (\Real \times \{0\}^2)$ and
$g_{i|B_i \cap (\Real \times \{0\}^2)} = g'_{|B_i \cap (\Real \times \{0\}^2)}$
\end{itemize}\end{defn}

Non-trivial connected-sums and prime knots are defined analogously.
A theorem of Schubert \cite{Sch} states that up to isotopy, 
every non-trivial $g$ can be written uniquely
up to a re-ordering of the terms, as a
connected-sum of prime knots $g=g_1 \# \cdots \# g_n$.

We review the Jaco-Shalen-Johannson decomposition of 3-manifolds \cite{JacoShalen}.
This is a standard decomposition of 3-manifolds along spheres and tori,
given by the connected-sum decomposition \cite{Kneser} followed by
the torus decomposition of the prime summands \cite{JacoShalen}
(see for example \cite{Hatcher3} or \cite{Neumann}). For us, all our
$3$-manifolds will be compact, and they are allowed to have a boundary.
For a more exhaustive treatment of JSJ-decompositions of knot and link
complements in $S^3$, see \cite{bjsj}.

\begin{figure}\label{fig9}
{
\psfrag{-1}[tl][tl][1][0]{$-1$}
\psfrag{i1}[tl][tl][0.7][0]{$-1+\frac{2}{2n+1}$}
\psfrag{i2}[tl][tl][0.7][0]{$-1+\frac{4}{2n+1}$}
\psfrag{i3}[tl][tl][0.7][0]{$-1+\frac{6}{2n+1}$}
\psfrag{i4}[tl][tl][0.7][0]{$-1+\frac{8}{2n+1}$}
\psfrag{im2}[tl][tl][0.7][0]{$-1+\frac{4n-2}{2n+1}$}
\psfrag{im1}[tl][tl][0.7][0]{$-1+\frac{4n}{2n+1}$}
\psfrag{1}[tl][tl][1][0]{$1$}
\psfrag{b1}[tl][tl][1][0]{$B_1$}
\psfrag{b2}[tl][tl][1][0]{$B_2$}
\psfrag{b3}[tl][tl][1][0]{$B_n$}
\psfrag{ddd}[tl][tl][1][0]{$\cdots$}
\psfrag{rx0}[tl][tl][1][0]{$\Real \times \{0\}^2$}
$$\includegraphics[width=13cm]{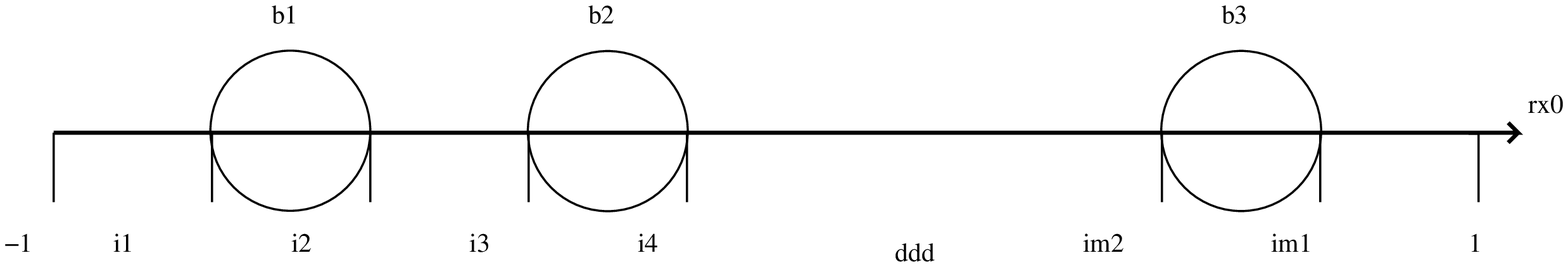}$$
\centerline{Figure 9}
}
\end{figure}

A $3$-manifold $M$ is a connected-sum $M=M_1 \# M_2$ if
surgery along an embedded $2$-sphere
produces manifolds $M_1$ and $M_2$, called the summands
of $M$. Provided
neither $M_1$ nor $M_2$ are $3$-spheres we say the connect-sum
is non-trivial. If a $3$-manifold $M$ is not $S^3$ and if all connect-sum
decompositions of $M$ are trivial, $M$ is called
prime. Kneser's Theorem \cite{Kneser} states that every compact, orientable
$3$-manifold is a connected-sum of
a unique collection of prime $3$-manifolds $M=M_1 \# M_2 \# \cdots \# M_n$,
where uniqueness is up to a re-ordering of the terms.

The torus decomposition of a prime
$3$-manifold $M$ consists of a minimal collection of embedded incompressible tori
$\sqcup_{i=1}^n T_i \subset M$ such that the complement
$M - \sqcup_{i=1}^{n} \nu T_i$ is a disjoint union of
atoroidal and Seifert-fibered manifolds, where $\nu T_i$ is an open tubular
neighborhood of $T_i\subset M$. A torus $T_i$ is incompressible if the induced
map $\pi_1 T_i \to \pi_1 M$ is injective.  A torus in a $3$-manifold
is peripheral if it is isotopic to a boundary torus. A $3$-manifold
is atoroidal if all incompressible tori are peripheral. 
The theorem of Jaco, Shalen and Johannson states that such a collection of
tori $\{T_1, T_2 \cdots, T_n\}$ always exists
and they are unique up to isotopy \cite{JacoShalen}.
Given an arbitrary prime $3$-manifold, there is an associated graph called
the JSJ-graph of $M$.  The vertices of the JSJ-graph are the components of
the manifold $M - \sqcup_{i=1}^{n} \nu T_i$.  The edges of the graph are
the tori $T_i$ for $i \in \{1,2,\cdots,n\}$.

Given a long knot $f \in \tK$, consider the compact
$3$-manifold $B-N'$ where $B\subset \Real^3$ is a closed $3$-ball containing
$\I\times D^2$, and $N'$ is the interior of the image of $f$.
 We will call $C=B-N'$ the knot complement. Define $T=\partial C$.
We review JSJ-splittings of knot complements.
Every sphere in $\Real^3$ bounds a $3$-ball by the Alexander-Schoenflies Theorem
(see for example \cite{Hatcher3}),
thus knot complements are prime $3$-manifolds, and the Jaco-Shalen-Johannson
decomposition of a knot complement is simply the torus decomposition.
The Generalized Jordan Curve
Theorem (see for example \cite{Pollack}) tells us a knot complement's
associated graph is a tree. The tree is rooted, as
only one component of $C-\sqcup_{i=1}^n \nu T_i$ contains $T$.
The component of $C-\sqcup_{i=1}^n \nu T_i$ containing $T$
will be called the root manifold of the JSJ-splitting.

\begin{defn} \label{pndef} Fix an embedding
 $b: \sqcup_{i \in \{1,2,\cdots,n\}} D^2 \to D^2$
such that $\partial D^2 \cap img(b)=\phi$. Let
$D^2_i$ denote the image of the $i$-th copy of $D^2$ under $b$.
Choose $b$ so that $D^2_i$ is the disc of radius $\frac{1}{2n+1}$ centered
around the point $(\frac{4i-2n-2}{2n+1},0)$.
Define $P_n$ to be $D^2 - int(\sqcup_{i=1}^n D^2_i)$. $P_n$ will be called the
$n$-times punctured disc.  $\partial D^2$ is the external boundary
and $\partial (img(b))$ the internal boundary of $P_n$ (see Figure 10).
\end{defn}

\begin{figure}\label{fig10}
{
\psfrag{-1}[tl][tl][1][0]{$-1$}
\psfrag{i1}[tl][tl][0.7][0]{$-1+\frac{2}{2n+1}$}
\psfrag{i2}[tl][tl][0.7][0]{$-1+\frac{4}{2n+1}$}
\psfrag{i3}[tl][tl][0.7][0]{$-1+\frac{6}{2n+1}$}
\psfrag{i4}[tl][tl][0.7][0]{$-1+\frac{8}{2n+1}$}
\psfrag{im2}[tl][tl][0.7][0]{$-1+\frac{4n-2}{2n+1}$}
\psfrag{im1}[tl][tl][0.7][0]{$-1+\frac{4n}{2n+1}$}
\psfrag{1}[tl][tl][1][0]{$1$}
\psfrag{d2}[tl][tl][1][0]{$D^2$}
\psfrag{b1}[tl][tl][1][0]{$D^2_1$}
\psfrag{b2}[tl][tl][1][0]{$D^2_2$}
\psfrag{b3}[tl][tl][1][0]{$D^2_n$}
\psfrag{ddd}[tl][tl][1][0]{$\cdots$}
\psfrag{rx0}[tl][tl][1][0]{$\Real \times \{0\}$}
$$\includegraphics[width=13cm]{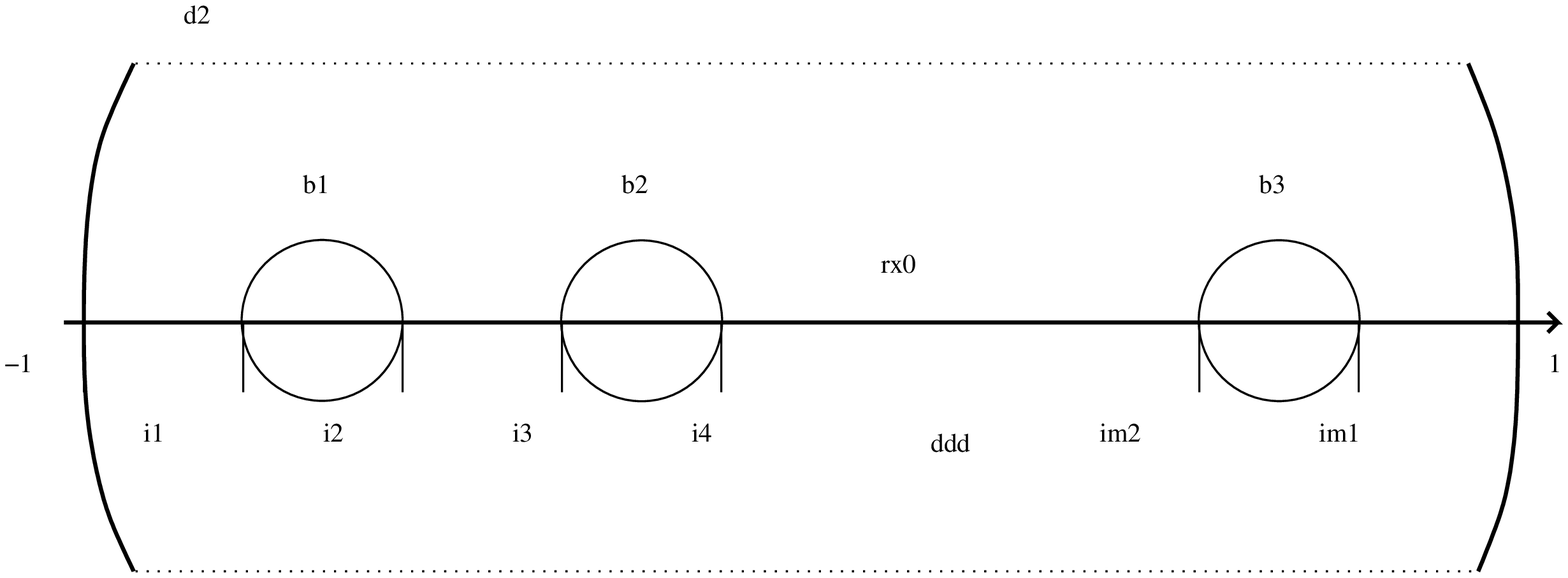}$$
\centerline{Figure 10}
}
\end{figure}

There are a few elementary facts that we will need about JSJ-splittings
of knot complements and diffeomorphism groups of $2$ and $3$-dimensional
manifolds. We assemble these facts in the following lemmas, all which are
widely `known' yet published proofs are elusive. A more detailed
study of JSJ-decompositions of knot and link complements in $S^3$ has
recently appeared \cite{bjsj} and could be used in place of several of
these lemmas. An essential reference for the following arguments is Hatcher's notes on $3$-dimensional manifolds \cite{Hatcher3}. 

\begin{lem}\label{sfl} If $M$ is sub-manifold of $S^3$ 
whose boundary consists of a non-empty collection of tori,
then either $M$ is a solid torus $S^1\times D^2$ or a component of the
complement of $M$ in $S^3$ is a solid torus.
\begin{pf}
Let $C=S^3 - int(M)$ be the complement. Since $\partial M$ consists of a disjoint union
of tori, every component of $\partial M$ contains an essential curve $\alpha$ which bounds a disc
$D$ in $S^3$.  Isotope $D$ so that it intersects $\partial M$ transversely
in essential curves. Then $\partial M \cap D \subset D$ consists of a finite collection
of circles, and these circles bound a nested collection of discs in $D$. Take an innermost
disc $D'$. If $D' \subset M$ then $M$ is a solid torus. If $D' \subset C$ then
the component of $C$ containing $D'$ is a solid torus.
\end{pf}
\end{lem}

\begin{lem}\label{sf2}
If a Seifert-fibred $3$-manifold is a component of a knot complement (in $S^3$) split
along its JSJ-decomposition, then it is diffeomorphic to one of the following:
\begin{itemize}
\item A solid torus (unknot complement).
\item The complement of a non-trivial torus knot. Such a
manifold is Seifert-fibred over a disc with two singular fibres.
\item $S^1 \times P_n$ for $n \geq 2$ (trivially fibred over a $n$-times punctured disc).
\item Fibred over an annulus with one singular fibre (The complement of
a regular and singular fibre in a Seifert-fibring of $S^3$).
\end{itemize}
\begin{pf}
Seifert-fibered manifolds that fiber over a non-orientable
surface do not embed in $S^3$ since a non-orientable, embedded closed
curve in the base lifts to a Klein bottle, which does not embed in $S^3$
by the Generalized Jordan Curve Theorem \cite{Pollack}. Similarly, a 
Seifert-fibered manifold that fibers over a surface of genus $g>0$ does
not embed in $S^3$ since the base manifold contains two curves that
intersect transversely at a point. If we lift one of these curves
to a torus in $S^3$, it must be non-separating. This again contradicts
the Generalized Jordan Curve Theorem.

Consider a Seifert-fibered manifold $M$ over an $n$-times punctured disc with
$n>0$ and with perhaps multiple singular fibers.
By Lemma \ref{sfl}, either $M$ is a solid torus or some component 
$Y$ of $\overline{S^3 \setminus M}$ is a solid torus. Consider 
the latter case.  There are two possibilities. 

\begin{enumerate}
\item The meridians of $Y$ are fibres of $M$. 
If there is a singular fibre in $M$, let $\beta$ be an embedded arc in the
base surface associated to the Seifert-fibring of $M$ which starts at the singular
point in the base and ends at the boundary component corresponding to $\partial Y$.
$\beta$ lifts to a $2$-dimensional CW-complex in $M$, and the endpoint of
$\beta$ lifts to a meridian of $Y$, thus it bounds a disc.  If we append this
disc to the lift of $\beta$, we get a CW-complex $X$ which consists of a $2$-disc
attached to a circle. The attaching map for the $2$-cell is multiplication by
$\beta$ where $\frac{\alpha}{\beta}$ is the slope associated to the singular fibre. 
The boundary of a regular neighbourhood of $X$ is a $2$-sphere, so we have 
decomposed $S^3$ into a connected sum $S^3 = L_{\frac{\gamma}{\beta}} \# Z$ 
where $L_{\frac{\gamma}{\beta}}$ is a lens space with $H_1 L_{\frac{\gamma}{\beta}} = \Zed_\beta$. Since $S^3$ is irreducible, $\beta=1$.
Thus $M\simeq S^1 \times P_{n-1}$ for some $n \geq 1$.  
\item The meridians of $Y$ are not fibres of $M$.
In this case, we can extend the Seifert fibring
of $M$ to a Seifert fibring of $M \cup Y$.  Either $M \cup Y = S^3$, or
$M \cup Y$ has boundary.
\begin{itemize}
\item If $M \cup Y = S^3$ then we know by the classification of Seifert fibrings of $S^3$ that any fibring of $S^3$ has at most two singular fibres.  If $M$ is the complement of
a regular fibre of a Seifert fibring of $S^3$, then $M$ is a torus knot complement.
Otherwise, $M$ is the complement of a singular fibre, meaning that $M$ is a solid torus.
\item If $M \cup Y$ has boundary, we can repeat the above argument. 
Either $M \cup Y$ is a solid torus, or a component of $S^3 \setminus \overline{M \cup Y}$ is a solid torus, so 
we obtain $M$ from the above manifolds by removing a Seifert fibre. By induction,
we obtain $M$ from either a Seifert fibring of a solid torus, or a Seifert 
fibring of $S^3$ by removing fibres. 
\end{itemize}
\end{enumerate}
\end{pf}
\end{lem}

\begin{figure}\label{fig11}
{
\psfrag{cabling}[tl][tl][0.9][0]{\hskip 6mm cabling}
\psfrag{jsj}[tl][tl][0.8][0]{JSJ}
\psfrag{tree}[tl][tl][0.8][0]{tree}
$$\includegraphics[width=13cm]{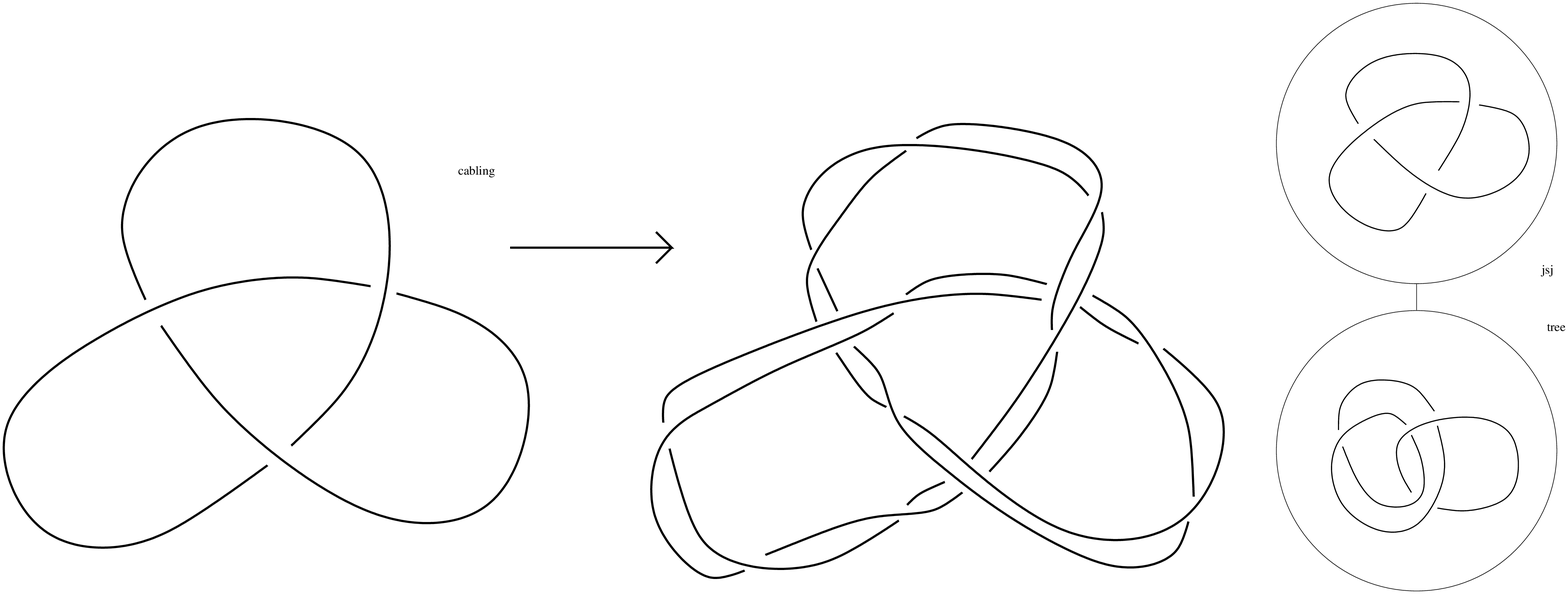}$$
\centerline{Figure 11}
}
\end{figure}

\begin{cor}\label{cablelem}
Given a long knot $f \in \tK$ with complement $C$, if the root manifold of the
JSJ-splitting of the knot complement is Seifert fibered
with one singular fiber, then $f$ is a cabling of another long knot.
Another way to say this is that the one-point compactification $\tilde f : S^1 \to S^3$ of 
$f_{|\Real \times \{0\}^2} : \Real \to \Real^3$ is an essential curve in
the boundary of a tubular neighborhood of some embedding $g : S^1 \to S^3$
(see Figure 11).
\end{cor}

Thurston \cite{Thurston} has proved that the non-Seifert-fibered manifolds in
the JSJ-splitting of a knot complement are finite-volume hyperbolic manifolds.
These hyperbolic manifolds can have arbitrarily many boundary components
\cite{bjsj}. 
Figure 12 demonstrates a hyperbolic satellite knot
(a knot such that the root manifold in the JSJ-decomposition is
hyperbolic) which contains the Borromean
rings complement in its JSJ-decomposition. 
In general, one can prove
that if the root manifold is a hyperbolic manifold with $n+1$ boundary components, 
then it is the complement of an $(n+1)$-component hyperbolic link in $S^3$ which contains
an $n$-component sublink which is the unlink. 

\begin{figure}\label{fig12}
{
\psfrag{n}[tl][tl][0.8][0]{A hyperbolic satellite operation}
\psfrag{jsj}[tl][tl][0.8][0]{JSJ}
\psfrag{tree}[tl][tl][0.8][0]{tree}
$$\includegraphics[width=10cm]{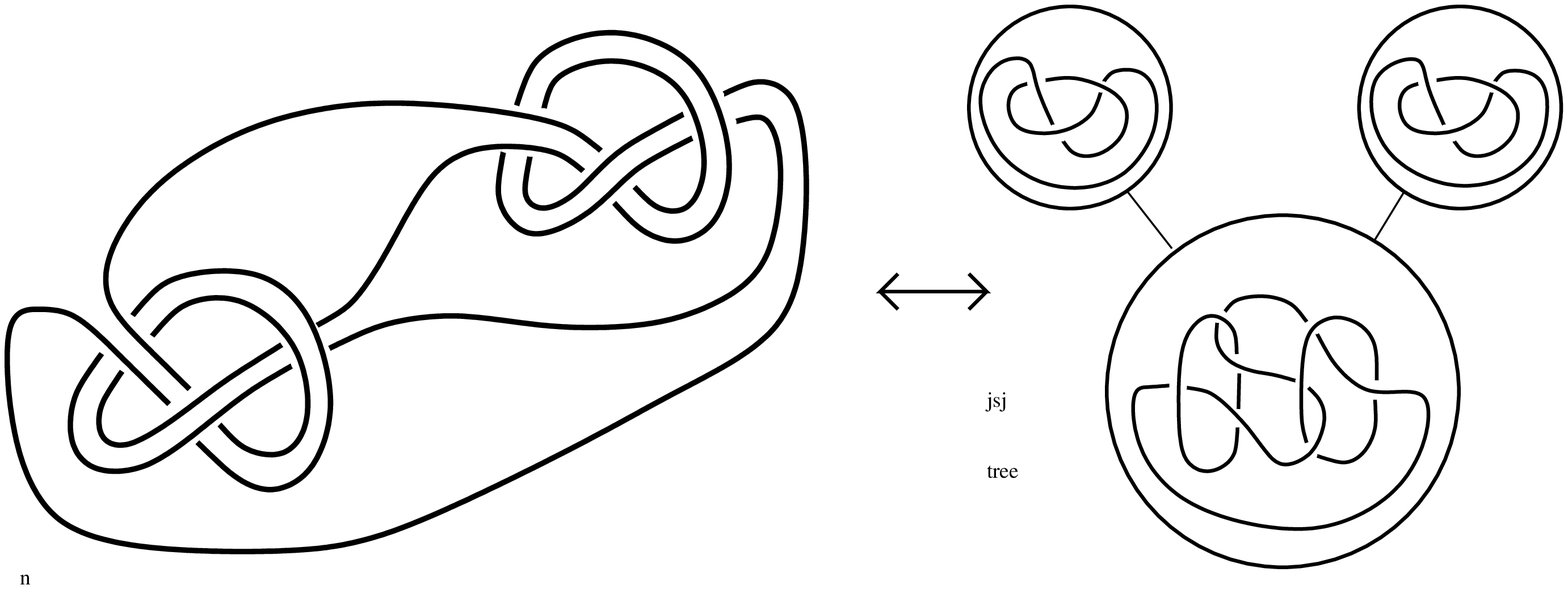}$$
\vskip 1mm
\centerline{Figure 12}
}
\end{figure}

\begin{lem}\label{pnlem}
A knot is a non-trivial connected-sum if and only if the root manifold
of the associated JSJ-tree
is diffeomorphic to $S^1 \times P_n$ for some $n \geq 2$. In this case,
$n$ is the number of prime summands of $f$.
\begin{pf} 
If $f \in \tK$ is a non-trivial connected-sum, let $n$ be
the number of prime summands of $f$, and isotope $f$ so that
$f_{|\Real\times \{0\}^2}$ satisfies Definition \ref{sumdef}.

Let $L \subset \Real^2$ be the closed disc of radius $\frac{1}{2}$ centered about
the origin. Let $N'=img(f_{|\Real \times L})$, $N=img(f)$
and define $C=B-int(N')$ where $B$ is a closed, convex ball neighborhood of
 $\I \times D^2$ in $\Real^3$.
Let $B_1, B_2, \cdots B_n$ be the closed $3$-balls from Definition
\ref{sumdef}, with $S_i = \partial B_i$ and $S_i$ intersecting
$img(f)$ in two discs for all $i \in \{1,2,\cdots,n\}$. Define
$C_i = B_i - int(N)$ and $T_i = \partial C_i$.

Let $\nu T_i$ be a small open tubular neighborhood of $T_i$, then
$C-\sqcup_{i=1}^n \nu T_i$ consists of $n+1$ components.
One component contains $T=\partial C$ and the other $n$ components are
the knot complements of the prime summands of $f$, $C_1,C_2,\cdots,C_n$.  
The component containing $T$ we will denote $V$.
$V$ is diffeomorphic to $S^1 \times P_n$.

By Dehn's Lemma the tori $\{T_i : i\in \{1,2,\cdots,n\}\}$ are incompressible in $C$.
If $\{T_{n+1},\cdots,T_{n+m}\}$ are the tori of the JSJ-decomposition for $\sqcup_{i=1}^n C_i$,
the collection
$\{T_1, T_2, \cdots, T_n, T_{n+1}, \cdots,  T_{n+m}\}$
is therefore the JSJ-decomposition of $C$.
Thus, $V \simeq S^1 \times P_n$ is the root manifold in the JSJ-tree associated to $C$.

To prove the converse, let $V$ be the root manifold of the JSJ-splitting of
$C$. Observe that $\partial V \simeq \partial (S^1 \times P_n)$ divides
$\Real^3$ into $n+2$ components, only one containing the knot.  Let
$T$ denote the boundary of the component which contains the knot.
By Lemma \ref{sf2} the fibers of $S^1 \times P_n$ are meridians of
the knot. Let $L_1,\cdots,L_n$
be properly embedded intervals in $P_n$ which cut $P_n$ into the union of
a disc with $n$ once-punctured discs.
Then $\sqcup_{i=1}^{n} (S^1 \times L_i)$ can be extended to $n$
disjoint, embedded $2$-spheres $S_i \subset \Real^3$ such that
$S_i \cap (S^1 \times P_n) = S^1 \times L_i$, and
$S_i \cap img(f_{|\Real\times \{0\}^2})$ consists of two points.
Thus we have decomposed the long knot $f$ into a connected-sum.
\end{pf}
\end{lem}

\begin{defn}\label{topm} In the above lemma, we call the tori $T_1, \cdots, T_n$ the
base level of the Jaco-Shalen-Johannson decomposition of the knot complement.
\end{defn}

Lemmas \ref{pnlem}, \ref{cablelem} and
Thurston's Hyperbolisation Theorem \cite{Thurston} gives us a canonical decomposition of
knots into simpler knots via cablings, connected-sums and hyperbolic satellite
operations commonly referred to as the satellite or splice decomposition of knots.  This is worked out in detail in \cite{bjsj}.

\begin{exmp}Figure 13 shows a knot with its JSJ-tori, and the associated JSJ-tree. In the standard terminology of knot theory, this knot would
be described as a connect-sum of three prime knots: the left-handed
trefoil, the figure-8 knot and the Whitehead double of the figure-8 
knot.  $V=S^1 \times P_3$ is the root manifold, $T_1,T_2,T_3$ are the base-level of the JSJ-decomposition of $C$, and $T_4$ is the 
remaining torus in the JSJ-decomposition of $C$. The leftmost 
summand is the trefoil knot.  The center summand is a figure-$8$
knot, whose complement is hyperbolic. The rightmost summand is the Whitehead double of the figure-$8$ knot, it's complement is $C_3$.  $C_3$ is 
the union of $C'_3$ (the Whitehead link complement, which
is hyperbolic) and $C''_3$ (a figure-$8$ knot
complement) where $\partial C'_3=T_3 \sqcup T_4$, and
$\partial C''_3=T_4$, $C''_3 \cap C'_3 = T_4$. The
interior of $C'_3$ is also a hyperbolic $3$-manifold of finite volume.
\end{exmp}

\begin{figure}\label{fig13}
{
\psfrag{k1}[tl][tl][0.8][0]{$T_1$}
\psfrag{k2}[tl][tl][0.8][0]{$T_2$}
\psfrag{k3}[tl][tl][0.8][0]{$T_3$}
\psfrag{k4}[tl][tl][0.8][0]{$T_4$}
$$\includegraphics[width=13cm]{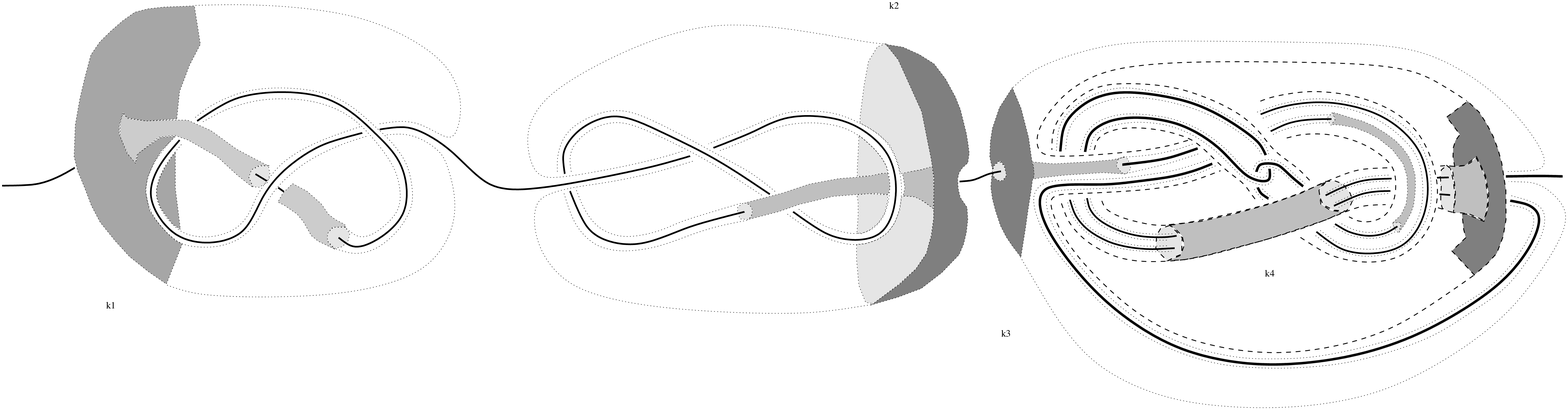}$$
}
{
\psfrag{s1d2}[tl][tl][1.2][0]{$V$}
\psfrag{t1}[tl][tl][1][0]{$T_1$}
\psfrag{t2}[tl][tl][1][0]{$T_2$}
\psfrag{t3}[tl][tl][1][0]{$T_3$}
\psfrag{tp}[tl][tl][1][0]{$T_4$}
\psfrag{whi}[tl][tl][1][0]{$C'_3$}
\psfrag{fig}[tl][tl][1][0]{$C_2$}
\psfrag{tre}[tl][tl][1][0]{$C_1$}
\psfrag{fig2}[tl][tl][1][0]{$C''_3$}
\psfrag{jsj}[tl][tl][1][0]{the JSJ}
\psfrag{tree}[tl][tl][1][0]{tree}
\psfrag{c3}[tl][tl][1][0]{$C'_3$}
$$\includegraphics[width=6cm]{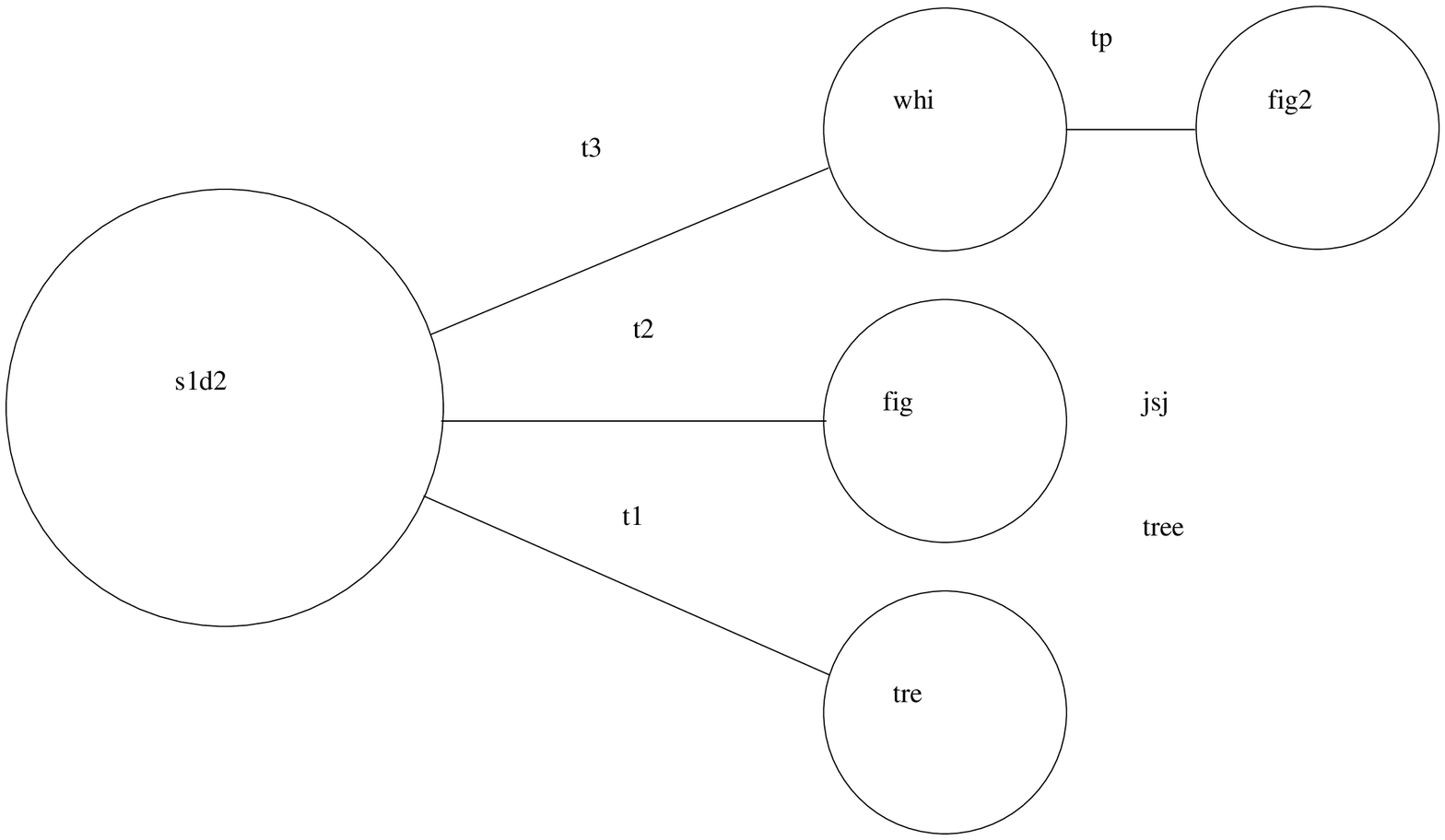} \hskip 1cm \includegraphics[width=3.4cm]{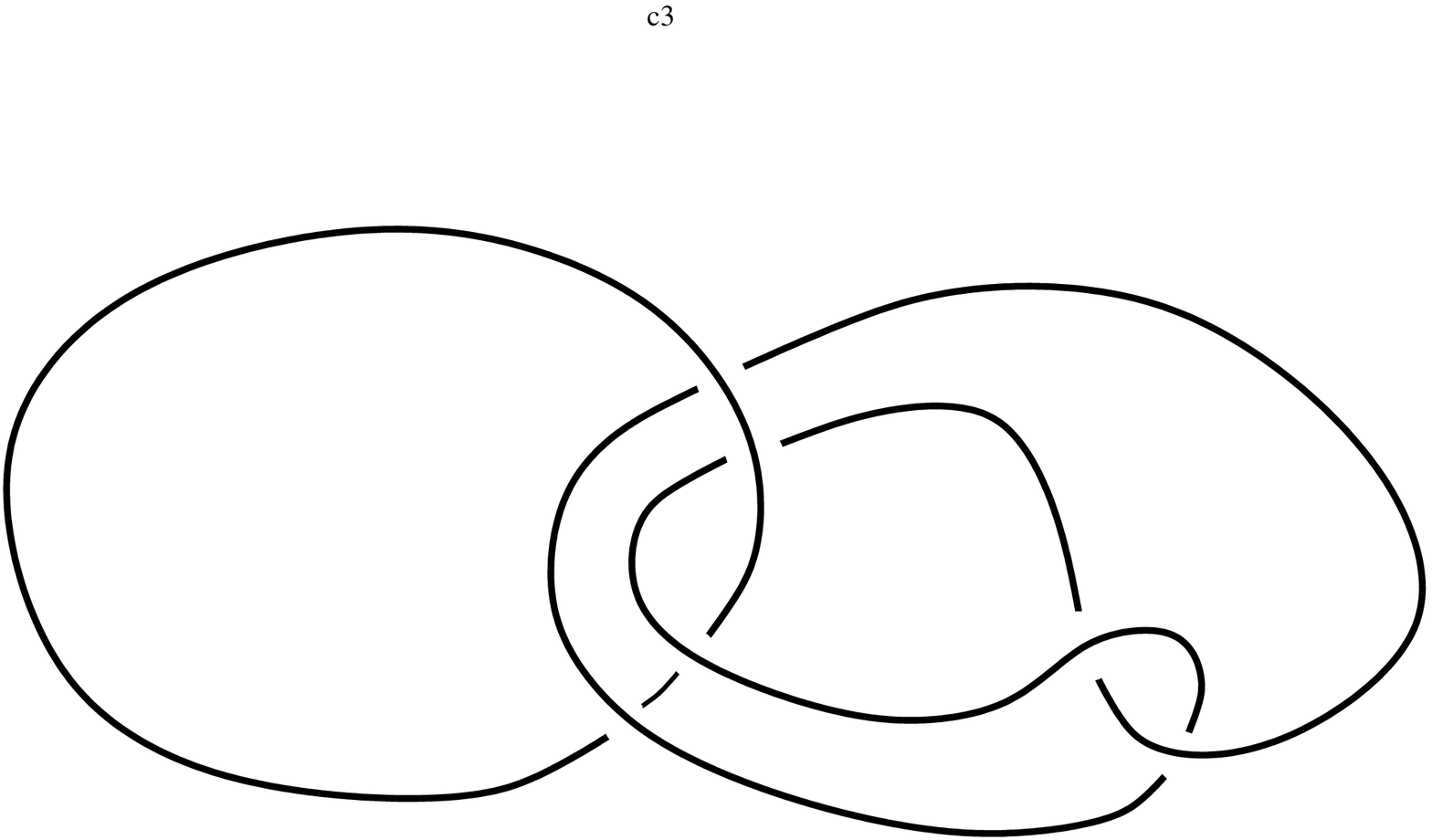}$$
\centerline{Figure 13}
}
\end{figure}

\begin{lem}\label{bdiff}
The component $\tK_f$ of $\tK$ containing the long knot $f$ is the classifying
space of $\Diff(C,T)$, the group of diffeomorphisms of the knot complement
which fix the boundary torus $T=\partial C$ point-wise. 
Moreover, $\tK_f$ is a $K(\pi,1)$.
\begin{pf}
 Let $B = \I\times D^2$ and let $\Diff(B)$ be the group of diffeomorphisms
 of $\Real^3$ with support contained in $B$.
 The map $\Diff(B) \to \tK_f$ defined by restriction to $img(f)$
 induces a fibration $$\Diff(C,T) \to \Diff(B) \to \tK_f$$
where $C=B-int(img(f))$, $T=\partial C$, and
$\Diff(C,T)$ is the group of diffeomorphisms of $C$ that fix $T$ point-wise.
Since $\Diff(B)$ is contractible \cite{Hatcher2}, $B\Diff(C,T) \simeq \tK_f$, where
$\tK_f$ is the component of $\tK$ containing $f$.
The fact that $\Diff(C,T)$ has contractible components is
due to Hatcher \cite{Hatcher1}.
\end{pf}
\end{lem}

In the above lemma, $BG=EG/G$ is the classifying space of a
topological group $G=\Diff(C,T)$ and $EG=\Diff(B)$.
Using Smale's Theorem $\Diff(D^2) \simeq \{*\}$ \cite{Smale}, 
an argument analogous to the above gives $\Cu_2(n)/S_n \simeq B\Diff(P_n)$ where $\Diff(P_n)$ is the group of diffeomorphisms of $P_n$ that
fix the external boundary of $P_n$ point-wise.

Let $\PDiff(P_n)$
denote the subgroup of $\Diff(P_n)$ consisting of
diffeomorphisms whose restrictions to $\partial P_n$ are isotopic to the identity
map $Id_{\partial P_n} : \partial P_n \to \partial P_n$. Then similarly,
by Smale's Theorem $\Cu_2(n) \simeq B\PDiff(P_n)$.

Let $\PFDiff(P_n)$ be the subgroup of $\PDiff(P_n)$ consisting of diffeomorphisms
whose restrictions to $\partial P_n$ are equal to the identity $Id_{\partial P_n}$.
$\pi_0 \Diff(P_n)$ is called the braid group on $n$-strands.
$\pi_0 \PDiff(P_n)$ is called the pure braid group on $n$-strands,
and $\pi_0 \PFDiff(P_n)$ is called the pure framed braid group on $n$
strands.  Observe that $\PFDiff(P_n)$ is homotopy equivalent to the subgroup $\PFDiff^+(P_n)$
of $\PFDiff(P_n)$ consisting of diffeomorphisms which restrict to the identity in an
$\epsilon$-neighborhood $N$ of the internal boundary of $P_n$.
  This follows from the fact that the space of
collar neighborhoods of $\partial P_n$ in $P_n$ is contractible.

\begin{defn}\label{windef} This definition will use the notation of Definition \ref{pndef} and
the previous paragraph.
Every diffeomorphism in $\PFDiff^+(P_n)$ can be canonically extended to a diffeomorphism
of the once-punctured disc $D^2 - int(D^2_i)$ simply by taking the union with
$Id_{D^2_j}$ for $j \neq i$.
Thus, for each $i \in \{1,2,\cdots,n\}$ there is a homomorphism
$w_i : \PFDiff^+(P_n) \to \pi_0 \PFDiff(S^1 \times \I) \simeq \Zed$ given by
the above extension together with an identification
$D^2 - int(D^2_i) \equiv S^1 \times \I$.  Here $\PFDiff(S^1 \times \I)$
denotes the group of boundary-fixing diffeomorphisms of $S^1 \times \I$. The generator
of $\pi_0 \PFDiff(S^1 \times \I)\simeq \Zed$ is a Dehn twist about a boundary-parallel
curve \cite{Gramain}.
Let $\dehn_n$ denote a free abelian subgroup of $\PFDiff^+(P_n)$ having rank $n$,
all whose elements have support in $N$, generated by Dehn twists about $n$ curves in $N$,
the $i$-th curve parallel to $\partial D^2_i$.
\end{defn}

\begin{lem}\label{pp}
There is an isomorphism of groups
$$ \pi_0 \PDiff(P_n) \times \Zed^n \simeq \pi_0 \PFDiff(P_n)$$
Moreover, the subgroups $\cap_{i=1}^n ker(w_i)$ and $\dehn_n$
satisfy:
\begin{itemize}
\item Inclusion $\cap_{i=1}^n ker(w_i) \to \PDiff(P_n)$ is a homotopy equivalence.
\item The elements of $\dehn_n$ and $\cap_{i=1}^n ker(w_i)$ commute with each other,
and $\dehn_n \cap \left(\cap_{i=1}^n ker(w_i) \right)$ is the trivial group.
\item The homomorphism $\cap_{i=1}^n ker(w_i) \times \dehn_n \to \PFDiff(P_n)$ is
a homotopy equivalence.
\end{itemize}
\begin{pf}
Take $\Diff(S^1)$ to be the group of orientation preserving diffeomorphisms of a circle,
and consider the fibration $\PFDiff(P_n) \to \PDiff(P_n) \to \prod_{i=1}^n \Diff(S^1)$
given by restriction to the internal boundary of $P_n$. This gives us the short exact
sequence

$$0 \to \prod_{i=1}^n \pi_1 \Diff(S^1) \to \pi_0 \PFDiff(P_n) \to \pi_0 \PDiff(P_n) \to 0$$

but $\prod_{i=1}^n \pi_1 \Diff(S^1) \simeq \Zed^n$, which is the subgroup
$\dehn_n \subset \pi_0 \PFDiff(P_n)$. The map
$\prod_{i=1}^n w_i : \pi_0 \PFDiff(P_n) \to \Zed^n \simeq \prod_{i=1}^n \pi_1 \Diff(S^1)$
is a splitting of the above short exact sequence. The kernel of $\prod_{i=1}^n w_i$ is
$\pi_0 \cap_{i=1}^n ker(w_i)$.  By definition, elements in $\cap_{i=1}^n ker(w_i)$ and $\dehn_n$
commute with each other, and so the result follows.
\end{pf}
\end{lem}

We will also need a mild variation on Lemma \ref{pp}.
 Let $*=(0,-1)$ be the base-point of $D^2$ and let $\gamma_i : [0,1] \to P_n$
for $i \in \{1,2,\cdots,n\}$ be
the affine-linear map starting at $*$ and ending at
$(\frac{4i-2n-2}{2n+1}, -\frac{1}{2n+1})$

\begin{figure}\label{fig14}
{\psfrag{d2}[tl][tl][1][0]{$$}
\psfrag{pn}[tl][tl][1][0]{$P_n$}
\psfrag{g1}[tl][tl][1][0]{$\gamma_1$}
\psfrag{g2}[tl][tl][1][0]{$\gamma_2$}
\psfrag{gn}[tl][tl][1][0]{$\gamma_n$}
\psfrag{star}[tl][tl][1][0]{$*=(0,-1)$}
\psfrag{rx0}[tl][tl][1][0]{$\Real \times \{0\}$}
\psfrag{cd}[tl][tl][1][0]{$\cdots$}
$$\includegraphics[width=5cm]{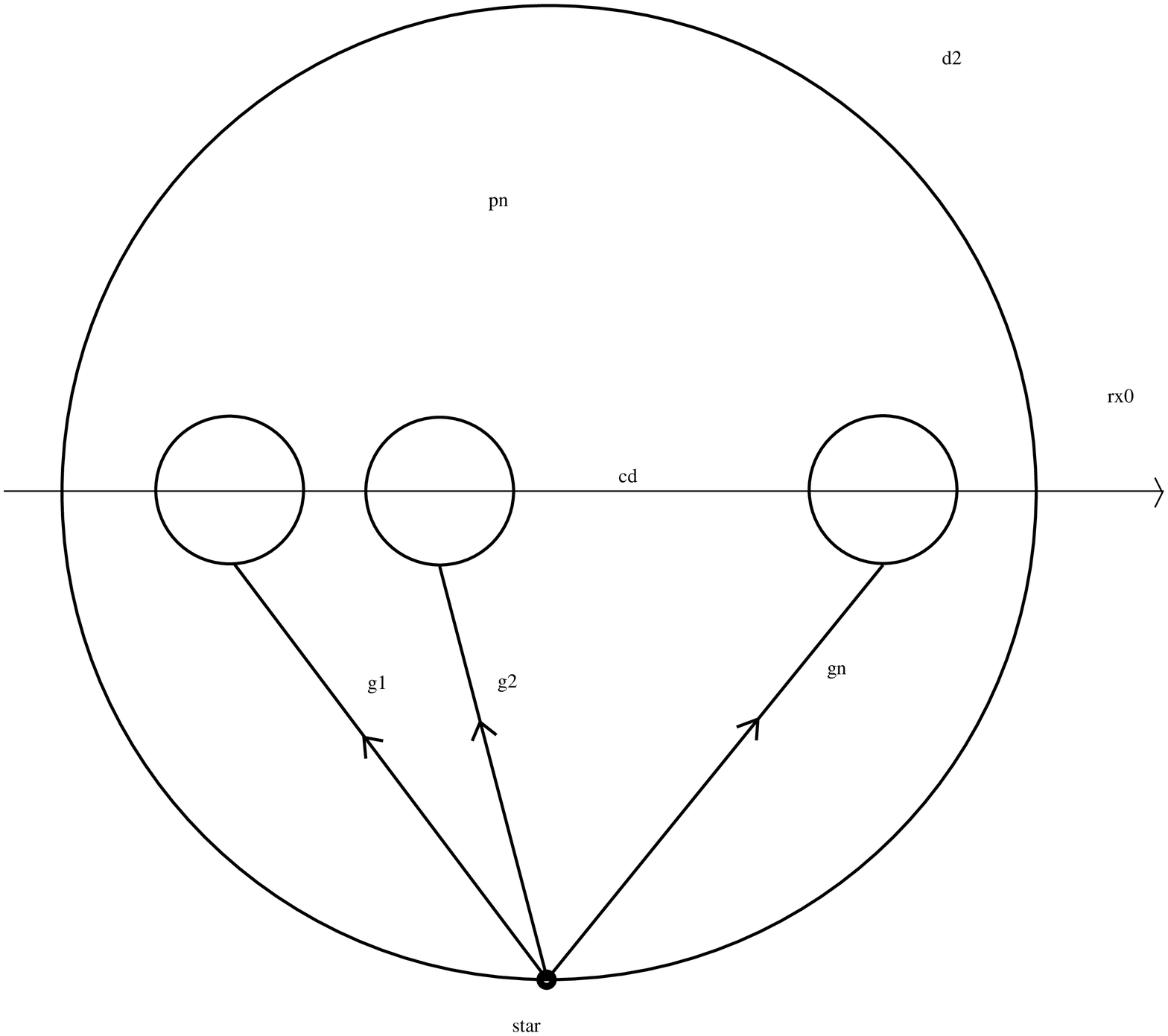}$$
\vskip 2mm
\centerline{Figure 14}
}
\end{figure}

\begin{defn} Define $\KDiff(P_n)$ to be $\cap_{i=1}^n ker(w_i)$. Define $\FDiff(P_n)$ to
be the subgroup of $\Diff(P_n)$ such  that each diffeomorphism $f \in \FDiff(P_n)$
\begin{itemize}
\item restricts to a diffeomorphism of $N$, ie: $f_{|N} : N \to N$.
\item the restriction of $f_{|N}$ to any connected component of $N$ is a translation
in the plane.
\end{itemize}
Observe, there is an epi-morphism $\FDiff(P_n) \to S_n \ltimes \Zed^n$ given by
$f \longmapsto (\sigma_f,\omega_1(f),\cdots,\omega_n(f))$ where
\begin{itemize}
\item $\sigma_f \in S_n$ is the permutation of $\{1,2,\cdots,n\}$ defined by  $\sigma_f(i)=j$
if $f(\partial D^2_i)=\partial D^2_j$.
\item $\omega_i(f) \in \Zed$ is the linking number of $\overline{\gamma_j}\cdot (f\circ \gamma_i)$ with
$D^2_j$ where $\sigma(i)=j$. Here $\overline{\gamma_j}(t)=\gamma_j (1-t)$ and concatenation is
by convention right-to-left, ie: if $\gamma, \eta : [0,1] \to X$ satisfy $\eta(1)=\gamma(0)$
then $\gamma \cdot \eta(t)=\eta(2t)$ for $0 \leq t \leq \frac{1}{2}$ and $\gamma\cdot\eta(t)=\gamma(2t-1)$ for
$\frac{1}{2} \leq t \leq 1$.
\item $S_n \ltimes \Zed^n$ is the semi-direct product of $S_n$ and $\Zed^n$ where $S_n$
acts on $\Zed^n$ by the regular representation ie:
$\sigma.(i_1,i_2,\cdots,i_n)=(i_{\sigma^{-1}(1)},i_{\sigma^{-1}(2)},\cdots,i_{\sigma^{-1}(n)})$
\end{itemize}
Call the above epi-morphism $W : \FDiff(P_n) \to S_n \ltimes \Zed^n$, and
define $\widetilde{\KDiff}(P_n) = W^{-1}(S_n \times \{0\}^n)$.
\end{defn}

\begin{lem}\label{wtkdiff} There is a fiber-homotopy equivalence
$$\xymatrix{
\PDiff(P_n) \ar[r] & \Diff(P_n) \ar[r] & S_n \\
\KDiff(P_n) \ar[r] \ar[u] & \widetilde{\KDiff}(P_n) \ar[r] \ar[u] & S_n \ar[u]
}$$
where all vertical arrows are inclusions.
\end{lem}

The above lemma follows immediately from Lemma \ref{pp}.

Abstractly there is a homotopy equivalence between
$B\KDiff(P_n)$ and $\Cu_2(n)$ given by the proof of Lemma \ref{bdiff}.
Since the properties of this homotopy equivalence will be important later, we
define it precisely here.
 
\begin{defn}\label{pd_def} Given $f \in \Diff(D^2)$, let $\zeta(f)=(L_1,L_2,\cdots,L_n) \in \Cu_2(n)$ 
be $n$ little $2$-cubes such that the center of $L_i$ is $f(\frac{4i-2n-2}{2n+1},0)$. 
For $\zeta(f)$ to be well-defined (and continuous) we need to choose the
the side lengths of $L_i$ equal to the minimum of these two numbers:
$\frac{1}{2n+1}$ and
the largest number $w$ so that the little cubes with centers
$f(\frac{4i-2n-2}{2n+1},0)$ with width and height equal to $w$ for 
$i \in \{1,2,\cdots,n\}$ have disjoint
interiors. Then $\phi : \Diff(D^2) \to \Cu_2(n)$ factors to a map
$B\KDiff(P_n) \to \Cu_2(n)$ which is a homotopy-equivalence.
\end{defn}

The definition below will use the conventions of Definition \ref{pndef},
in particular we will call $S^1 \times \partial D^2 \subset S^1 \times P_n$ the external
boundary of $S^1 \times P_n$, and $\partial (S^1 \times P_n) - S^1 \times \partial D^2$
the internal boundary of $S^1 \times P_n$.

\begin{defn}\label{s1xpndiff}
Let $\eta_i : S^1 \to \partial D^2_i$ be a clockwise parametrization of $\partial D^2_i$
starting and ending at $\gamma_i(1)$. Notice that $\lambda_i = \overline{\gamma_i}\eta_i\gamma_i$
for $i \in \{1,2,\cdots,n\}$ are generators for $\pi_1 P_n$.
Let $\{*\}\times \lambda_i$ and $S^1 \times \{*\}$
denote generators of $\pi_1 (S^1 \times \partial D^2_i)$.
Let $\Diff(S^1 \times P_n)$ be the group of diffeomorphisms of $S^1 \times P_n$
whose restriction to the external boundary are equal to the identity
$Id_{S^1 \times \partial D^2}$
and whose restriction to the internal boundary $S^1 \times \partial(img(b))$
sends $\{1\}\times \eta_i$ to a curve isotopic to $\{1\} \times \eta_{\sigma(i)}$
for all $i \in \{1,2,\cdots, n\}$  where
 $\sigma : \{1,2,\cdots,n\} \to \{1,2,\cdots,n\}$ is a permutation of $\{1,2,\cdots,n\}$.
Let $\PDiff(S^1 \times P_n)$ denote the group of diffeomorphisms of
$S^1 \times P_n$ whose restrictions to the internal boundary are isotopic
to the identity and whose restrictions to the external boundary are
equal to the identity $Id_{S^1 \times \partial D^2}$.
Similarly, define $\PFDiff(S^1 \times P_n)$ to be the group of diffeomorphisms
of $S^1 \times P_n$ which restrict to the identity $Id_{S^1\times \partial P_n}$.
Let $\KDiff(S^1 \times P_n)$ be the subgroup of $\PFDiff(S^1 \times P_n)$
consisting of diffeomorphisms having the
form $Id_{S^1} \times f$ where $f \in \KDiff(P_n)$, and let $\widetilde{\KDiff}(S^1 \times P_n)$
denote the subgroup of $\Diff(S^1 \times P_n)$ consisting of diffeomorphisms of the form
$f=Id_{S^1} \times g$ for $g \in \widetilde{\KDiff}(P_n)$.
\end{defn}

\begin{lem}\label{diffs1pn}
 There is a fiber-homotopy equivalence
$$\xymatrix{
\PDiff(S^1 \times P_n) \ar[r] & \Diff(S^1 \times P_n) \ar[r] & S_n \\
\KDiff(S^1 \times P_n) \ar[r] \ar[u] & \widetilde{\KDiff}(S^1 \times P_n) \ar[r] \ar[u] & S_n \ar[u] }$$

where all vertical arrows are inclusions (and homotopy equivalences).
\begin{pf}
We consider $S^1 \times P_n$ to be a Seifert fibered manifold.
Hatcher \cite{Hatcher1} proves that
the full group of diffeomorphism of $S^1 \times P_n$ is homotopy
equivalent to the fiber-preserving subgroup.
Let $G$ denote the fiber-preserving subgroup of $\PDiff(S^1 \times P_n)$.
Thus, the inclusion $G \to \PDiff(S^1 \times P_n)$ is a homotopy equivalence.
Since the group of orientation preserving diffeomorphisms of $S^1$
is homotopy equivalent to $SO_2$,
$G$ is homotopy equivalent to the subgroup $G' \subset G$
of fibrewise-linear diffeomorphisms of $S^1 \times P_n$.
Since every diffeomorphism in $\PDiff(S^1 \times P_n)$
restricts to a diffeomorphism of $\partial (S^1 \times P_n)$ which
is isotopic to the identity, $G'$ is homotopy equivalent to the
subgroup of diffeomorphisms of the form $Id_{S^1} \times f$ where
$f \in \PDiff(S^1 \times P_n)$. The key consideration in the
above argument is whether
or not $f$ could be a Dehn twist along a vertical annulus.
By Lemma \ref{pp}, $\PDiff(P_n)$ is homotopy equivalent to $\KDiff(P_n)$.
The remaining results follow from Lemmas \ref{wtkdiff} and \ref{pp}.
\end{pf}
\end{lem}

As a historical note, some of Hatcher's results on diffeomorphism groups
of Haken manifolds were independently discovered by Ivanov \cite{Ivanov1, Ivanov2}.

The following lemma is used to simplify
the proof of Theorem \ref{freeness}.  It is a standard variation of
a construction of Borel \cite{Borel} (chapter IV, \S 3).

\begin{lem}\label{cohenlemma}
If $G$ is a topological group with $H$ a closed normal subgroup such that
$ G/H $ is a finite group, then there
exists a canonical normal, finite-sheeted covering space
$$ G/H \to BH \to BG$$
where the map $BH \to BG$ is given by
the projection $EG/H \to EG/G$ where we make the identification $BH=EG/H$.
\end{lem}

First, we sketch the proof of Theorem \ref{freeness}.
The fact that the map $\sqcup_{n=0}^\infty \kappa_n$ induces a bijection
$$\sqcup_{n \in \{ 0,1,2,3,\cdots \}}
 \pi_0 \left( \left( \Cu_2(n) \times \Prime^n\right)/S_n \right) \to \pi_0 \tK$$
is due to Schubert \cite{Sch}.  His theorem states that every long
knot decomposes uniquely into a connected-sum of prime knots, up
to a re-ordering of the terms.
Since the map $\sqcup_{n=0}^\infty \kappa_n$ is bijective on components,
we need only to verify
that it is a homotopy equivalence when restricted to any single
connected component.
By Lemma \ref{bdiff}, the components of both the domain and range
are $K(\pi,1)$'s.  So we have reduced the theorem to checking that the
induced map is an isomorphism of fundamental groups for every component.
The inspiration for the proof of this is the fibration below,
 which we call the little cubes fibration.

$$ S_n \to \Cu_2(n) \times \Prime^n \to (\Cu_2(n) \times \Prime^n)/S_n $$

Let $f \in \tK$ with $f=f_1\# f_2 \# \cdots \# f_n$, where $(f_1,\cdots,f_n)\in \Prime^n$
are the prime summands of $f$. Let $\tK_f$ denote
the component of $\tK$ containing $f$, similarly define $\tK_{f_i}$.
Thus the above fibration, when restricted to the connected component
$\Cu_2(n) \times \prod_{i=1}^n \tK_{f_i}$ of $\Cu_2(n) \times \Prime^n$, has the form:

$$\Sigma_f \to \Cu_2(n) \times \prod_{i=1}^n \tK_{f_i} \to
(\Cu_2(n) \times \prod_{i=1}^n \tK_{f_i})/\Sigma_f$$

where $\Sigma_f \subset S_n$ is the subgroup which preserves the partition
$\sim $ of $\{1,2,\cdots,n\}$
with $i \sim j \Leftrightarrow \tK_{f_i} = \tK_{f_j}$.

By Lemma \ref{bdiff} the little cubes fibration gives the short exact
sequence below.

$$ 0 \to \pi_1 \Cu_2(n) \times \prod_{i=1}^n \pi_1 \tK_{f_i} \to
         \pi_1 ((\Cu_2(n) \times \prod_{i=1}^n \tK_{f_i})/\Sigma_f) \to \Sigma_f \to 0$$

 $\pi_1 \tK_f \simeq \pi_0 \Diff(C,T)$ by Lemma
\ref{bdiff}.  So the idea of the proof is to find an analogous fibration
for $\tK_f$.  So we are looking for an epimorphism
$\pi_0 \Diff(C,T) \to \Sigma_f$.

Since the tori in the JSJ-splitting of $C$ are unique up to isotopy, define
a permutation $\sigma_g : \{1,2,\cdots,n\} \to \{1,2,\cdots,n\}$ by the
condition that $\sigma_g(i)=j$ if $g(T_i)$ is isotopic to $T_j$ where
$T_1,T_2, \cdots,T_n$ are the base-level of the JSJ-decomposition of $C$.
This is well-defined since $g$ fixes $T = \partial C$ and the JSJ-decomposition
is unique up to isotopy.
The homomorphism
$\sigma : \pi_0 \Diff(C,T) \to S_n$ is onto $\Sigma_f$ since
two long knots $f_i$ and $f_j$ are isotopic if and only if $C_i$ and $C_j$
admit orientation preserving diffeomorphisms which also preserve the
(oriented) meridians of $C_i$ and $C_j$.

The kernel of $\sigma$
one would expect to be the mapping class group of diffeomorphisms
of $C$ which do not permute the base-level
of the JSJ-splitting of $C$. Such a diffeomorphism $g$, when
restricted to $V \simeq S^1 \times P_n$ can
isotoped to be in $\KDiff(S^1 \times P_n)$. Thus $g$ 
restricts to diffeomorphisms $g_{|C_i} \in \Diff(C_i,T_i)$
for all $i \in \{1,2,\cdots, n\}$, leading us to expect the kernel of $\sigma$
is $\pi_0 \PDiff(P_n) \times \prod_{i=1}^n \pi_0 \Diff(C_i,T_i)$.
By Lemma \ref{bdiff} $\pi_0 \Diff(P_n) \simeq \pi_1 \Cu_2(n)$ and
$\pi_0 \Diff(C_i,T_i)\simeq \pi_1 \tK_{f_i}$
where $f_i$ denotes the $i$-th summand of $f$.  So we have constructed
a SES
 $$ 0 \to \pi_1 \Cu_2(n)\times \prod_{i=1}^n \pi_1 \tK_{f_i} \to \pi_1 \tK_f \to \Sigma_f \to 0$$
which is the analogue of the SES coming from the little cubes fibration.

In the argument below, we rigorously redo the above sketch
at the space-level. We construct a fibration of
 diffeomorphism groups whose long exact sequence is the SES given above.
We then use Lemma \ref{cohenlemma} to convert this fibration of
 diffeomorphism groups into a fibration which describes $\tK_f$, and this
 we will show is equivalent to the little cubes fibration.


\begin{pf} (of Theorem \ref{freeness})

We will show that $\sqcup_{n=0}^\infty \kappa_n$ is a homotopy
equivalence, component by component.
Let $f \in \tK$ be a knot specifying a connected
component $\tK_f$ of $\tK$.

In the case of the unknot $f = Id_{\Real\times D^2}$, we know from
the proof of the Smale conjecture \cite{Hatcher2} that the component of $\tK$ containing
$f$ is contractible.  $\Cu_2(0) \times \Prime^0$ is a point thus
the map $\kappa_0$ is a homotopy equivalence between these
two components.

If $f$ is a prime knot, $n=1$ and the  little cubes
fibration $S_1 \to \Cu_2(1)\times \Prime^1 \to (\Cu_2(1)\times \Prime^1)/S_1$
is trivial, thus $\tK_f$ is a component of $\Prime$.
In this case, our map $\kappa_1 : \Cu_2(1) \times \Prime \to \tK$ is mapping from $\Cu_2(1) \times \Prime$
to $\tK$. Since $\Cu_2(1)$ is contractible and our action
satisfies the identity axiom, $\kappa_1$ is homotopic to the
composite of the projection map $\Cu_2(1) \times \Prime \to \Prime$
with the inclusion map $\Prime \to \tK$, and so $\kappa_1$ is a homotopy
equivalence between $(\Cu_2(1) \times \Prime)/S_1$ and $\Prime$.

Consider the case of a composite knot $f=f_1 \# f_2 \# \cdots \# f_n \in \tK$
for $n \geq 2$ with $f_i$ prime for all $i \in \{1,2,\cdots,n\}$.
Let $C = B - N'$ denote the knot complement, as in Lemma \ref{pnlem}.
Let $T=\partial C$, let $V \simeq S^1 \times P_n$
denote the root manifold of the associated tree to the JSJ-decomposition
of $C$ and let $T_1,\cdots,T_n$ denote base-level of the JSJ-decomposition of
$C$ (see Lemma \ref{pnlem}, Definition \ref{topm}). Similarly, let
$V \simeq S^1 \times P_n$, $B_i$ and $C_i$ for $i \in \{1,2,\cdots,n\}$
be as in Lemma \ref{pnlem}.
Let $\Diff(C,T)$ be the group of diffeomorphisms of $C$ that fix $T$
point-wise. Let $\Diff^V(C,T)$ denote the subgroup of $\Diff(C,T)$
consisting of diffeomorphisms which
restrict to diffeomorphisms of $V$.  Let $\PDiff^V(C,T)$ denote the subgroup
of $\Diff^V(C,T)$ consisting of diffeomorphisms
whose restrictions to $\partial V$ are isotopic
to $Id_{\partial V}$.
Let $\Emb(\sqcup_{i=1}^n T_i,C)$  denote the space of embeddings of
$\sqcup_{i=1}^n T_i$ in $C$.
If we restrict a diffeomorphism in $\Diff(C,T)$ to $\sqcup_{i=1}^n T_i$
and mod-out by the parametrization of the individual tori, we get a fibration
(which is not necessarily onto)

$$\PDiff^V(C,T) \to \Diff(C,T) \to \Emb(\sqcup_{i=1}^n T_i,C)/\prod_{i=1}^n \Diff(T_i)$$

Since $T_i$ is incompressible in $C$, this fibration is
mapping to embeddings
which are also incompressible. The tori $\sqcup_{i=1}^n T_i$
are part of the JSJ-splitting of
$C$, and the JSJ-splitting is unique up to isotopy. This
means that a diffeomorphism
in $\Diff(C,T)$ must send $T_i$ to another torus in the
 JSJ-splitting (up to isotopy), but more
importantly that torus must be in the base-level of the JSJ-splitting since
the diffeomorphism is required to preserve $T$.

A component of $\Emb(\sqcup_{i=1}^n T_i,C)/\prod_{i=1}^n \Diff(T_i)$
is an isotopy class of $n$ embedded, labeled tori.  Provided the tori
are incompressible, such a component must be contractible \cite{Hatcher1}.
Consider the union $X$ of all the components of
$\Emb(\sqcup_{i=1}^n T_i,C)/\prod_{i=1}^n \Diff(T_i)$
which correspond to embeddings whose image are the base-level of the JSJ-splitting of $C$.
$X$ must have the homotopy type of the symmetric group $S_n$.
Consider $S_n$ to be the subspace
$S_n \equiv \Diff(\sqcup_{i=1}^n T_i)/\prod_{i=1}^n \Diff(T_i)
 \simeq X \subset \Emb(\sqcup_{i=1}^n T_i,C)/\prod_{i=1}^n \Diff(T_i)$.

The above argument proves that there is a fiber-homotopy
equivalence, where all the vertical arrows are given by inclusion.

$$\xymatrix{\PDiff^V(C,T) \ar[r] & \Diff(C,T) \ar[r] & X \\
            \PDiff^V(C,T) \ar[r] \ar[u] & \Diff^V(C,T) \ar[r] \ar[u] &  S_n \ar[u] }$$

Typically it is demanded that fibrations are onto. Since the long knot $f$ is a connected-sum,
and some of the summands $\{f_i : i \in \{1,2,\cdots,n\}\}$ may be repeated, define
the equivalence relation $\sim$ on $\{1,2,\cdots, n\}$ by $i \sim j \Leftrightarrow 
f_i$ is isotopic to $f_j$. Let $\Sigma_f \subset S_n$ be the partition-preserving
subgroup of $S_n$. Thus the above fibration is onto $\Sigma_f \subset S_n$.

Since every diffeomorphism $g \in \PDiff^V(C,T)$ restricts
to a diffeomorphism of $V$, consider the restriction to $V \simeq S^1 \times P_n$.
Since the $g$ extends to a diffeomorphism of $\Real^3$,
$g_{|V} : V \to V$ must preserve (up to isotopy) the longitudes and meridians
of each $T_i$. To be precise, a meridian of $T_i$ is an oriented closed essential curve in
$T_i$ which bounds a disc in $\Real^3 - int(C_i)$. The orientation of the meridian
is chosen so that the
linking number of the meridian with the knot is $+1$. A longitude in $T_i$ is an essential
oriented curve in $T_i$ which bounds a Seifert surface in $C_i$. The
orientation of the curve is chosen to agree with the orientation of $f_i$.

Thus, if we identify $V$ with $S^1 \times P_n$ in a way that sends
knot meridians to fibers of $S^1 \times P_n$ and the longitude of
$f_i$  to
$\{1\}\times \eta_i\subset S^1 \times P_n$ for all $i \in \{1,2,\cdots,n\}$
 then (by a slight abuse of notation) 
 $g_{|S^1 \times P_n} \in \PDiff(S^1 \times P_n)$.

Define $\KDiff^V(C,T) \subset \PDiff^V(C,T)$ 
   and $\widetilde{\KDiff}^V(C,T) \subset \Diff^V(C,T)$
to be the subgroups such that each diffeomorphism $g$ restricts to a diffeomorphism
of $V \equiv S^1 \times P_n$,
$g_{|S^1 \times P_n} \in \KDiff(S^1 \times P_n)$ and 
$g_{|S^1 \times P_n} \in \widetilde{\KDiff}(S^1 \times P_n)$ respectively.
By Lemma \ref{diffs1pn}, the vertical inclusion maps in the diagram below
give a fiber-homotopy equivalence

$$\xymatrix{\PDiff^V(C,T) \ar[r] & \Diff^V(C,T) \ar[r] & \Sigma_f \\
            \KDiff^V(C,T) \ar[u]^\simeq
	    \ar[r] \ar[u]^\simeq & \widetilde{\KDiff}^V(C,T) \ar[r] 
	    \ar[u]^\simeq & \Sigma_f \ar[u]^\simeq }$$

Analogously to Lemma \ref{pp}, the
inclusion $\KDiff(S^1 \times P_n) \times \prod_{i=1}^n \Diff(C_i,T_i) \to \KDiff^V(C,T)$
is a homotopy equivalence.


If we apply Lemma \ref{cohenlemma} to the above fibration, we get the
normal covering space

$$ \xymatrix{\Sigma_f \ar[r] &
             B\KDiff(S^1 \times P_n) \times \prod_{i=1}^n B\Diff(C_i,T_i) \ar[r] \ar[d]_\simeq &
	     B\widetilde{\KDiff}^V(C,T) \ar[d]_\simeq \\
             & \Cu_2(n)\times \prod_{i=1}^n \tK_{f_i}  & \tK_f}$$

where the two vertical homotopy equivalences come from Lemma \ref{bdiff}  and the
identification $\KDiff(S^1 \times P_n) \equiv \KDiff(P_n)$

Consider $\Cu_2(n) \times \prod_{i=1}^n \tK_{f_i}$ as a $\Sigma_f$-space,
where the $\Sigma_f$ action is simply the restriction of the diagonal $S_n$ action
$S_n \times (\Cu_2(n) \times \tK^n) \to \Cu_2(n) \times \tK^n$ to 
$\Sigma_f \times (\Cu_2(n) \times \prod_{i=1}^n \tK_{f_i}) \to \Cu_2(n) \times \prod_{i=1}^n \tK_{f_i}$.
  By design, the homotopy equivalence
$B\KDiff^V(S^1 \times P_n) \times \prod_{i=1}^n B\Diff(C_i,T_i) \to
 \Cu_2(n) \times \prod_{i=1}^n \tK_{f_i}$
is $\Sigma_f$-equivariant (see Definition \ref{pd_def}).

Thus we know abstractly that there exists a homotopy equivalence between
$(\Cu_2(n)\times \prod_{i=1}^n \tK_{f_i})/\Sigma_f$ and $\tK_f$.
To finish the proof, we show
 $\kappa_n :(\Cu_2(n)\times \prod_{i=1}^n \tK_{f_i})/\Sigma_f \to \tK_f$
is such a homotopy equivalence.
Since both the domain and range of $\kappa_n$ are $K(\pi,1)$'s, it suffices to show that the
diagram below commutes.

$$ \xymatrix{\pi_1 B\KDiff(S^1 \times P_n) \times \prod_{i=1}^n \pi_1 B\Diff(C_i,T_i) \ar[r] \ar[d]_\simeq &
	     \pi_1 B\widetilde{\KDiff}^V(C,T) \ar[d]_\simeq \\
             \pi_1 \Cu_2(n)\times \prod_{i=1}^n \pi_1 \tK_{f_i} \ar[r]^-{\pi_1 \kappa_n} & \pi_1 \tK_f}$$

Fix $i \in \{1,2,\cdots,n\}$ and $\phi \in \pi_0 \Diff(C_i,T_i)$. Consider
$\phi$ to be an element of
$\pi_1 B\KDiff(S^1 \times P_n) \times \prod_{i=1}^n \pi_1 B\Diff(C_i,T_i)$
by the standard inclusion.
If one chases $\phi$ along the clockwise route around the diagram to $\pi_1 \tK_f$,
one is simply converting $\phi$ into an element $\overline \phi \in \pi_1 \tK_f$
using Lemma \ref{bdiff}.  This means that one is applying an isotopy to
the $i$-th knot summand $f_i$ of $f$, and the isotopy has support in
$B_i$ (see Lemma \ref{pnlem}). If one chases $\phi$ along the counter-clockwise
route around the diagram, one converts $\phi$ into a loop in $\pi_1 \tK_{f_i}$ using
Lemma \ref{bdiff}, then the little cubes construction is applied to this loop creating
a second loop $\tilde \phi \in \pi_1 \tK_f$. The
loop produced via the little cubes construction $\tilde \phi$ is the same
loop in $\pi_1 \tK_f$ as $\overline \phi$ since the little cubes and other
knot summands remain fixed through the isotopy, keeping the support of the isotopy
in $B_i$.

Given $\theta \in \pi_0 \KDiff(S^1 \times P_n)$
consider it as an element of
$\pi_1 B\KDiff(S^1 \times P_n) \times \prod_{i=1}^n \pi_1 B\Diff(C_i,T_i)$ by the
standard inclusion.
We will chase $\theta$ around the diagram. This chase is a little more involved than
the previous one, as it involves the little cubes action on $\tK$ in a non-trivial manner.

Our strategy for the proof is to chase $\theta$ around the diagram in a counter-clockwise
manner to get an element in $\pi_0 \widetilde{\KDiff}^V(C,T)$. We denote this
diffeomorphism by $C_\theta$.
We need to show that $C_\theta$ is the identity on $\sqcup_{i=1}^n C_i$ and
when restricted to $V$, $C_{\theta | V}\equiv \theta$ under our identification
$V \equiv S^1 \times P_n$.
We will do this via an explicit computation. First, notice that we can simplify
the problem.  $\theta$ determines a loop $\widetilde{\theta} \in \pi_1 \Cu_2(n)$
which in turn defines an isotopy $\kappa_n(\widetilde{\theta},f_1,f_2,\cdots,f_n)$
of $f$, which by Lemma \ref{bdiff} determines the diffeomorphism $C_\theta$ of $C$.
 Recall how
$C_\theta$ is constructed. Given an isotopy $F_\theta: [0,1] \times B \to B$ such that
\begin{itemize}
\item $F_\theta(0,x)=x$ for all $x \in B$
\item $F_\theta(t,x)=x$ for all $x \in T = \partial B$ and $t \in [0,1]$
\item $F_\theta(t,x)=\kappa_n(\widetilde{\theta}(t),f_1,f_2,\cdots,f_n)(x)$ for all
      $(t,x) \in [0,1] \times B$.
\end{itemize}
Then $C_\theta(x)=F_\theta(1,x)$ for $x \in C$.

Define $T_\theta : B \to B$ by $T_\theta(x)=F_\theta(1,x)$ for $x \in B$.
$\pi_0 \KDiff(S^1 \times P_n) \simeq \pi_0 \KDiff(P_n)$ is the pure braid group
which can be in turn thought of as a subgroup of the full braid group,
 $\pi_0 \widetilde{\KDiff}(S^1 \times P_n) \simeq
 \pi_0 \widetilde{\KDiff}(P_n) \simeq \pi_0 \Diff(P_n)$.
In $\pi_0 \Diff(P_n)$ every element can be written as a product of Artin generators
$\{\sigma_i : i \in \{1,2,\cdots,n-1\}\}$ (see for example \cite{Birman}), these are
the half Dehn twists about curves bounding the $i$-th and $(i+1)$-st punctures of $P_n$.
Let
$\theta=\alpha_j \circ \alpha_{j-1} \circ \cdots \circ \alpha_1$
where $\alpha_i \in \Diff(P_n)$ are either Artin generators or their
inverses, thus
$T_\theta = T_{\alpha_j} \circ T_{\alpha_{j-1}} \circ \cdots \circ T_{\alpha_1}$.
This in principle reduces our problem to studying
$T_{\sigma_i}$ for $i \in \{1,2,\cdots,n-1\}$.

\begin{figure}\label{fig15}
{
\psfrag{1}[tl][tl][1][0]{$1$}
\psfrag{n}[tl][tl][1][0]{$n$}
\psfrag{i}[tl][tl][1][0]{$i$}
\psfrag{i+1}[tl][tl][1][0]{$i+1$}
\psfrag{s}[tl][tl][1][0]{$\sigma_i$}
$$\includegraphics[width=10cm]{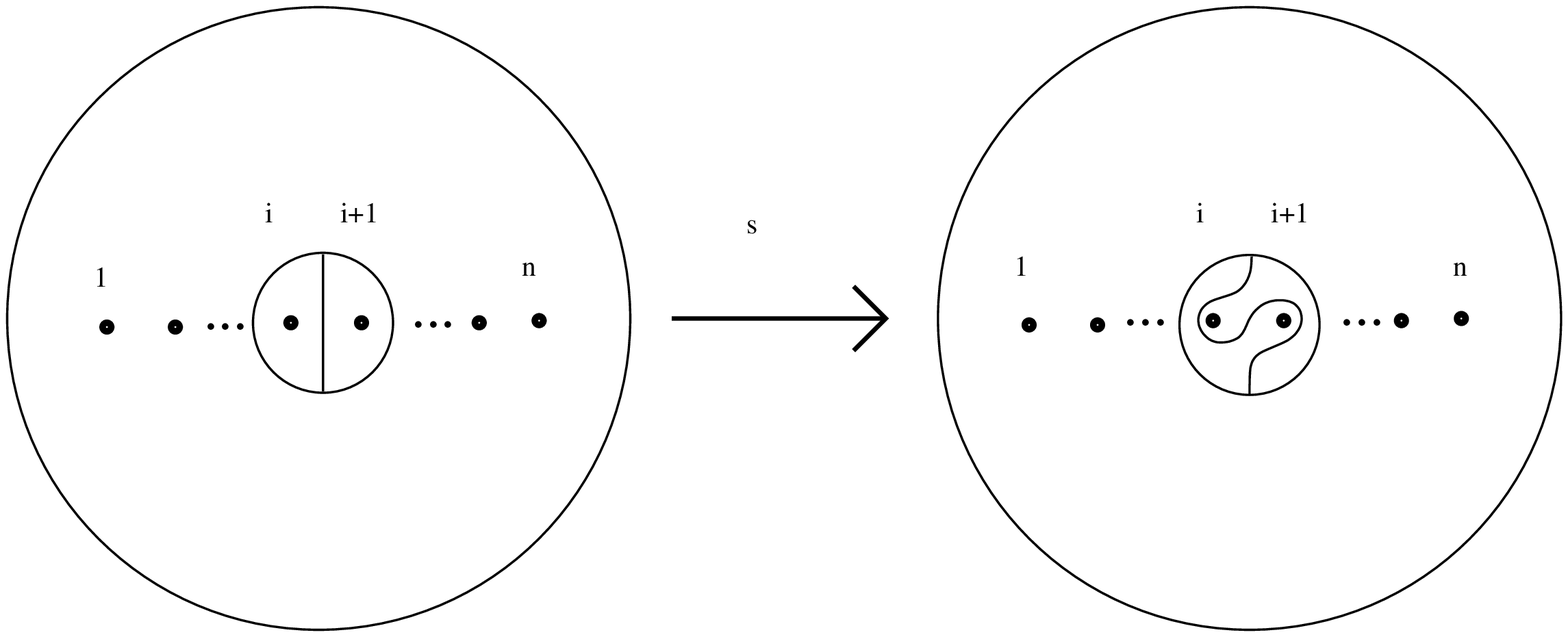}$$
\centerline{Figure 15}
}
\end{figure}
\begin{figure}\label{fig16}
{
\psfrag{1}[tl][tl][0.9][0]{$1$}
\psfrag{2}[tl][tl][0.9][0]{$i$}
\psfrag{3}[tl][tl][0.9][0]{$i+1$}
\psfrag{4}[tl][tl][0.9][0]{$n$}
\psfrag{cd}[tl][tl][0.9][0]{$\cdots$}
\psfrag{s}[tl][tl][1][0]{$\widetilde{\sigma_{i}}$}
$$\includegraphics[width=10cm]{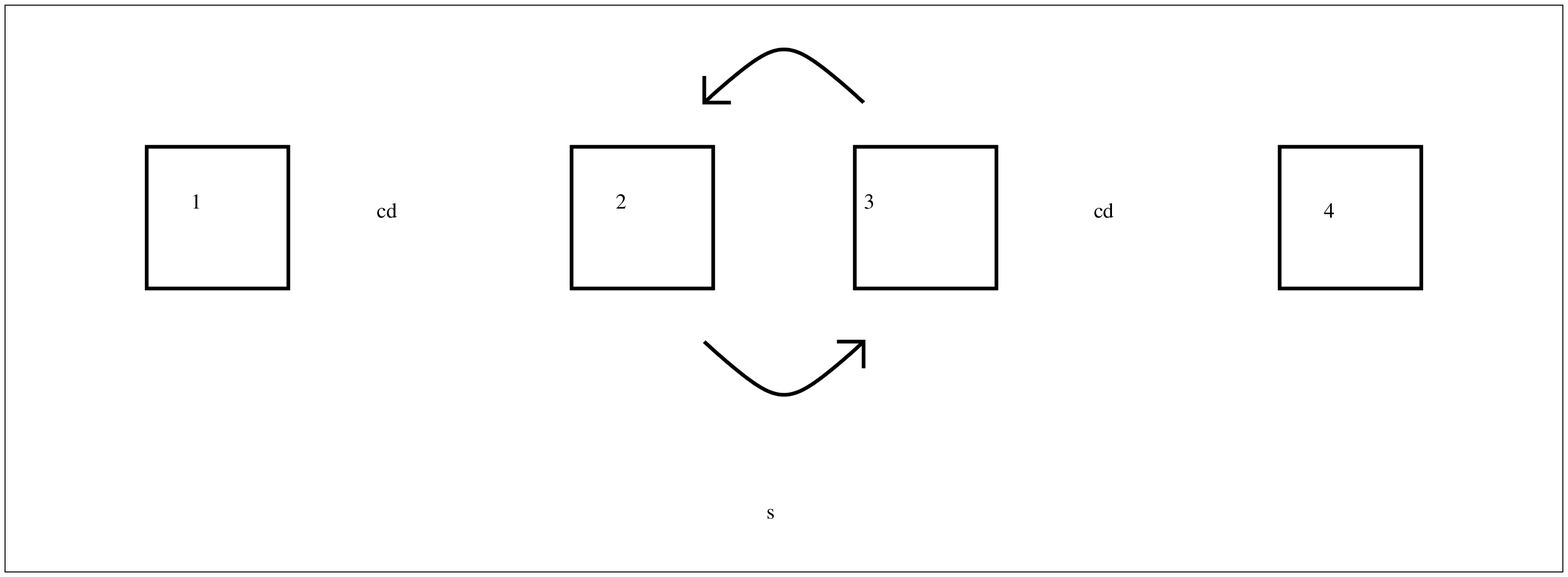}$$
\centerline{Figure 16}
}
\end{figure}

  By the definition of $\kappa_n$, $T_{\sigma_i}$
is the identity on the balls $B_k$ for $k \notin \{i,i+1\}$, and $T_{\sigma_i}$
permutes the two balls $B_i$ and $B_{i+1}$, acting by translation.  Thus
$T_\theta$ must restrict to be the identity on $\sqcup_{i=1}^n C_i$.

Let $*=(0,-1,0) \in \partial B$ be the base-point of $B$.  Let $\xi_i : [0,1] \to B$ be the
unique affine-linear function so that $\xi_i(0)=*$ and
$\xi_i(1)=(\frac{4i-2n-2}{2n+1},-\frac{1}{2n+1},0)\in \partial B_i$. Let
$p_i : S^1 \to C_i$ be a longitude of $C_i$ starting and ending at $\xi_i(1)$.
Since $T_{\sigma_i}$ acts by translation on the balls $\{B_i : i \in \{1,2,\cdots,n\}\}$,
for all $k \in \{0,1,2,\cdots,j\}$ define the $i$-th longitude $p_i^k$ of
$(\alpha_k \circ \alpha_{k-1} \circ \cdots \circ \alpha_1)(C)$
to be the restriction of
$\sqcup_{s=1}^n \left(\alpha_k \circ \alpha_{k-1} \circ \cdots
 \circ \alpha_1 \circ p_s\right) :
\sqcup_{s=1}^n S^1 \to B$
to $(\sqcup_{s=1}^n \alpha_k \circ \alpha_{k-1} \circ \cdots \circ \alpha_1 \circ p_s)^{-1}(B_i)$.
Define $l_i = \overline{\xi_i} \cdot p_i \cdot \xi_i$ and similarly
$l_i^k=\overline{\xi_i} \cdot p_i^k \cdot \xi_i$, so $l_i^0=l_i=l_i^j$ for all
$i \in \{1,2,\cdots,n\}$.

$\pi_1 \left( (\alpha_k \circ \alpha_{k-1} \circ \cdots \circ \alpha_1)(C)\right)$
therefore has a natural identification with
$\Zed \times (*_{i=1}^n \Zed)$ which has presentation
$\langle m,l^k_1,l^k_2,\cdots,l^k_n : [m,l^k_1],[m,l^k_2],\cdots,[m,l^k_n]\rangle$.
Here $m$ is a knot meridian, or equivalently a fiber of the Seifert
fibering of the base-manifold of the JSJ-splitting of $C$.

\begin{figure}\label{fig17}
{
\psfrag{l1}[tl][tl][0.9][0]{$l_1$}
\psfrag{l2}[tl][tl][0.9][0]{$l_2$}
\psfrag{ln}[tl][tl][0.9][0]{$l_n$}
\psfrag{b1}[tl][tl][0.9][0]{$B_1$}
\psfrag{b2}[tl][tl][0.9][0]{$B_2$}
\psfrag{bn}[tl][tl][0.9][0]{$B_n$}
\psfrag{star}[tl][tl][0.9][0]{$*$}
\psfrag{dot}[tl][tl][0.9][0]{$\cdots$}
$$\includegraphics[width=10cm]{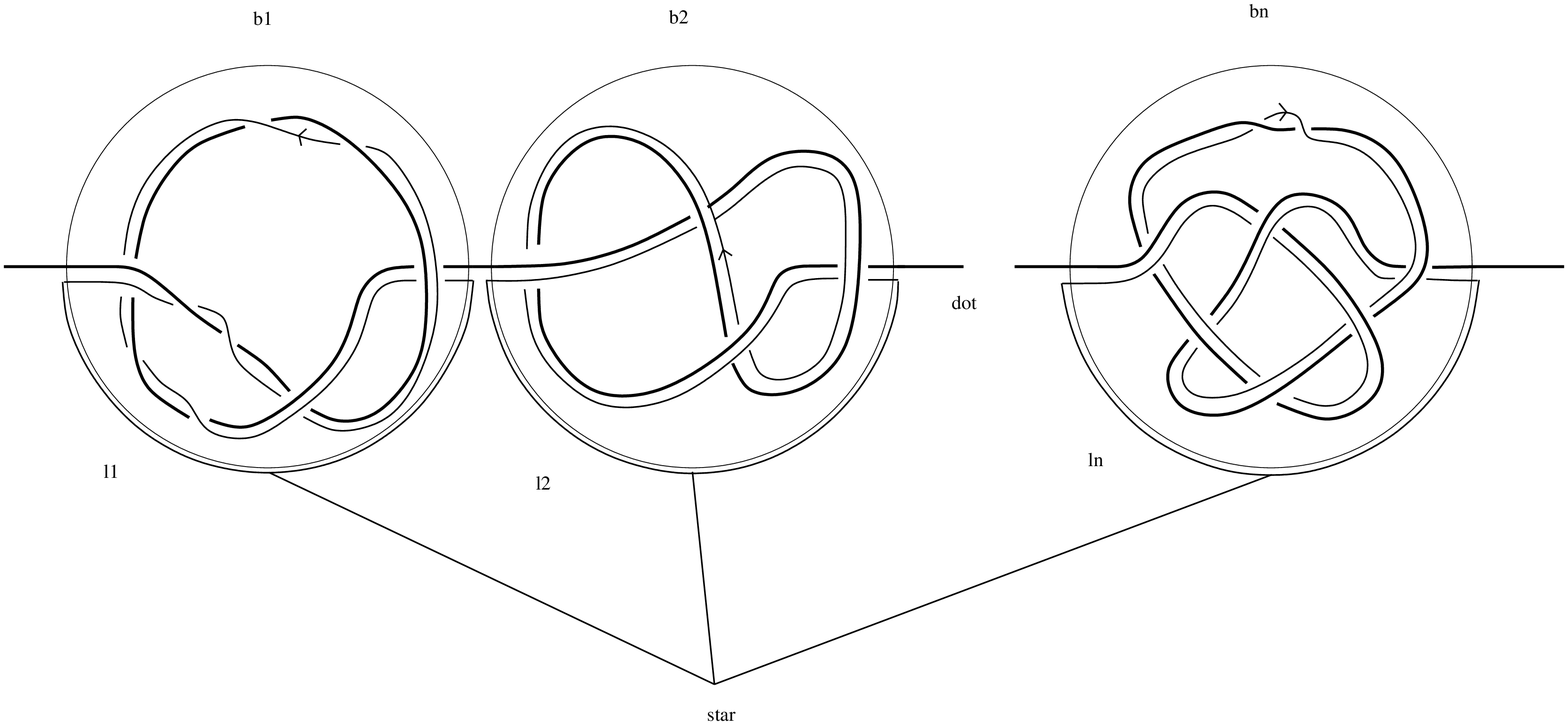}$$
\centerline{Figure 17}
}
\end{figure}

Call the above identification
$\phi_k : \pi_1 \left( (\alpha_k \circ \alpha_{k-1} \circ \cdots \circ \alpha_1)(C)\right)
 \to \Zed \times (*_{i=1}^n \Zed)$. $\phi_k$ determines a diffeomorphism
$\widetilde{\phi_k} : (\alpha_k \circ \alpha_{k-1} \circ \cdots \circ \alpha_1)(C)
 \to S^1 \times P_n$
defined by the condition that $\widetilde{\phi_k}(l^k_i)=\{1\} \times \lambda_i$,
$\widetilde{\phi_k}(m)=S^1 \times \{*\}$.

Recall the Dehn-Nielsen theorem \cite{Nielsen} (see \cite{Zie} for a modern proof).
It states that the map $\pi_0 \Diff(P_n) \to \Aut(\pi_1 P_n)$ is injective. We compute the
induced automorphism on $\Zed \times (*_{i=1}^n \Zed)$ given by the composite
$\widetilde{\phi_{k+1}} \circ T_{\alpha_{k+1}} \circ \widetilde{\phi_k^{-1}}$.
Without loss of generality, assume $\alpha_{k+1} = \sigma_q$ for some $q \in \{1,2,\cdots,n-1\}$,
therefore $\kappa_n (\widetilde{\alpha_{k+1}},f_1,f_2,\cdots,f_n)$ represents an isotopy
which pulls the knot summand in the ball $B_{q+1}$ through the knot summand in the ball
$B_q$. Therefore,
$\pi_1 \left(\widetilde{\phi_{k+1}} \circ T_{\alpha_{q}} \circ \widetilde{\phi_k^{-1}}\right)$
fixes $m$
and fixes $\lambda_i$ unless $i \in \{q,q+1\}$,
in which case
$(\widetilde{\phi_{k+1}} \circ T_{\alpha_{q}} \circ \widetilde{\phi_k^{-1}})(\lambda_{q+1})=
 \lambda_{q+1}\lambda_{q}\lambda_{q+1}^{-1}$
and 
$(\widetilde{\phi_{k+1}} \circ T_{\alpha_{q}} \circ \widetilde{\phi_k^{-1}})(\lambda_{q+1})= 
\lambda_q$.

Thus, via our identifications, 
$C_\theta \in \KDiff^V(C,T)$ induces the same automorphism of
$\pi_1 V \equiv \pi_1 (S^1 \times P_n)$ as
does $\theta \in \KDiff(S^1 \times P_n)$, which proves the theorem.
\end{pf}

\begin{cor} There is a little $2$-cubes equivariant homotopy equivalence
$$\EK{1,D^2} \simeq \Cu_2 (\Prime \sqcup \{*\}) \times \Omega^2\Bbb CP^\infty$$
where $\Bbb CP^\infty=BS^1=B^2\Zed$.
\end{cor}

\section{Where from here?}\label{endsec}

There are several directions one could go from here.  One direction
would be to ask, what is the homotopy type of the full space $\K$?
By Theorem \ref{freeness} this is equivalent to asking what is the
homotopy type of $\Prime$ but Theorem \ref{freeness} can be used
to refine this question further. 

 Starting with
the unknot, one can produce new knots by: using hyperbolic satellite
operations, cablings, or taking the connected-sum of knots.  If these procedures
are iterated, one produces all knots \cite{Thurston, JacoShalen, bjsj}. Theorem \ref{freeness}
tells us the homotopy-type of a component corresponding to a
knot which is a connect-sum.
If $f \sim f_1 \# \cdots \# f_n$ is the prime decomposition of $f$, then
$\K_f \simeq (\Cu_2(n) \times_{S_n} \prod_{i=1}^n \K_{f_i})$.
To complete
our understanding of $\K$ all we need to understand is:
 
\begin{enumerate}
\item How the homotopy type of $\K_f$ is related to
the homotopy type of $\K_g$ if $f$ is a cabling of $g$.
\item If $f$ is
obtained from knots $\{f_i : i \in \{1,2,\cdots,n\}\}$ via a 
hyperbolic satellite operation, how is
the homotopy type of $\K_f$ related to $\K_{f_i}$ for $i \in \{1,2,\cdots, n\}$.
\end{enumerate}

Hatcher has answered question 1.

\begin{thm} \label{htc} (Hatcher) \cite{Hatcher4}
If a knot $f$ is a cabling of a knot $g$ then
 $\K_f \simeq S^1 \times \K_g$
\end{thm}

More recently, a solution to question 2 has appeared in \cite{topknot}.
Roughly, if a knot $f$ is obtained from knots
$\{f_i : i \in \{1,2,\cdots,n\}\}$ by a hyperbolic satellite
operation then there is a fibration

$$\prod_{i=1}^n \K_{f_i} \to \K_f \to S^1 \times S^1$$

and the monodromy of this fibration depends on both the knots
$f_i$, their symmetry properties, and the symmetry properties of 
the hyperbolic manifold that is the root of the JSJ-tree of $f$. 
For brevity, we skip the full statement of the result.  A
key theorem of Sakuma's is used to compute the monodromy of
this fibration -- allowing us to show the fibration is split
at the base, thus the fundamental group of any component
of $\K$ is an iterated semi-direct product of finite-index
subgroups of braid groups.

More generally, one could ask, what is the homotopy type of
other spaces of knots?  

Perhaps the next simplest case is the space of embeddings of a circle in
a sphere $\Emb(S^1,S^n)$.  As is shown in \cite{BudCoh}, there is a homotopy
equivalence $\Emb(S^1,S^n) \simeq \Emb(\Real,\Real^n) \times_{SO_{n-1}} SO_{n+1}$.
Thus, if one knows the homotopy type of $\Emb(\Real,\Real^n)$ as an $SO_{n-1}$-space,
one knows the homotopy type of $\Emb(S^1,S^n)$. The homotopy-type of 
$\K$ as an $SO_2$-space is determined in \cite{topknot}.

Another interesting question is `what is the homotopy type of the space of
closed, connected, $1$-dimensional submanifolds of $S^n$'? This space is naturally 
homeomorphic to $\Emb(S^1,S^n)/\Diff(S^1)$ and has been studied recently by Hatcher \cite{Hatcher4} in the $n=3$ case. Studying the homotopy type of these spaces appears to have
more complications due to the delicate extension problems involved.
An interesting point of Hatcher's work is that one needs
to know the answer to the Linearization Conjecture in order to understand
even the homotopy type of the component of a knot as simple as a hyperbolic
knot. One could go further and ask, what is the homotopy-type of the
double-coset space $SO_{n+1}\backslash \Emb(S^1,S^n)/\Diff(S^1)$? This
is a particularly delicate problem as the action of $SO_{n+1}\times \Diff(S^1)$
on $\Emb(S^1,S^n)$ is not free. A nice example of the kinds of problems that
can arrise is the paper of Kodama and Michor \cite{komi}, where they prove that the
figure-8 component of $\Imm(S^1,\Real^2)/\Diff(S^1)$ has the homotopy-type of
$\Bbb CP^\infty$.

It would be very interesting to know more about the homotopy-type of
the embedding spaces $\Emb(\Real^j,\Real^n)$ or $\Emb(S^j,S^n)$.
Unfortunately the techniques of this paper are of limited use since it is
still unknown whether or not a smooth embedded $3$-sphere
in $\Real^4$ bounds a smooth ball \cite{Kirby}, and very little is known about
the homotopy type of $\Diff(D^4)$ other than Morlet's `Comparison' Theorem
\cite{Mor, BL, KS}.

There are however some results known in dimension $4$.  Sinha and Scannell
prove have computed many rational homotopy-groups of the long
knot space $\Emb(\Real,\Real^4)$ and the corresponding framed long
knot space $\EK{1,D^3}$, showing non-triviality
in dimensions $\{2,4,5,6\}$. The fibration $\Diff(D^4) \to \EK{1,D^3}$ 
has a fiber which is homotopy equivalent
to $\Diff(S^2 \times D^2)$ (diffeomorphisms fixing the boundary). The
homotopy LES of this fibration splits into short exact sequences
$0 \to \pi_{i+1} \EK{1,D^3} \to \pi_i \Diff(S^2\times D^2) \to \pi_i \Diff(D^4) \to 0$.
We can deduce from this that $\pi_i \Diff(S^2 \times D^2)$ has non-torsion
elements for $i\in \{1,3,4,5\}$. By Theorem \ref{littlecthm} we know
that $\Diff(S^2 \times D^2) \simeq \EK{2,S^2}$ is a $3$-fold loop space.
Three-dimensional instincts might lead one to
suspect that the inclusion $\Omega^2 SO_3 \subset \Diff(S^2 \times D^2)$
is a homotopy equivalence,  where $\Omega^2 SO_3$ is thought of as
the subgroup of
fiber-preserving (fibrewise-linear) diffeomorphisms of $S^2 \times D^2$.
These instincts would be wrong! 
We have just seen that although the inclusion 
$\Omega^2 SO_3 \to \Diff(S^2 \times D^2)$ admits a $3$-fold de-looping,
it can not be a homotopy equivalence since the homotopy groups of the
domain and range are not the same.

A possible application of the Sinha, Scannell result would be the study of
`spun' knots.  Given $f \in \pi_i \Emb(\Real,\Real^n)$ one constructs
a smooth embedding $S^{i+1} \to \Real^{n+i}$ by `spinning' $f$ about
an $(n-1)$-dimensional linear subspace of $\Real^{n+i}$ 
(this is a slight generalization
of Litherland's notion of deform twist-spun knots \cite{Lith}, 
see Figure 18). 
In the spirit of Markov's Theorem \cite{Birman}, it would seem natural to conjecture that
for some co-dimensions the `spinning map' 
$$\pi_i \Emb(\Real,\Real^n) \to \pi_0 \Emb(S^{i+1},\Real^{n+i})$$
\noindent is an isomorphism.  
\begin{figure}\label{fig18}
{
\psfrag{spun}[tl][tl][1][0]{}
\psfrag{r4}[tl][tl][1][0]{$\Real^{i+n}$}
\psfrag{r3}[tl][tl][1][0]{$\Real^{n-1}$}
\psfrag{t}[tl][tl][1][0]{}
$$\includegraphics[width=8cm]{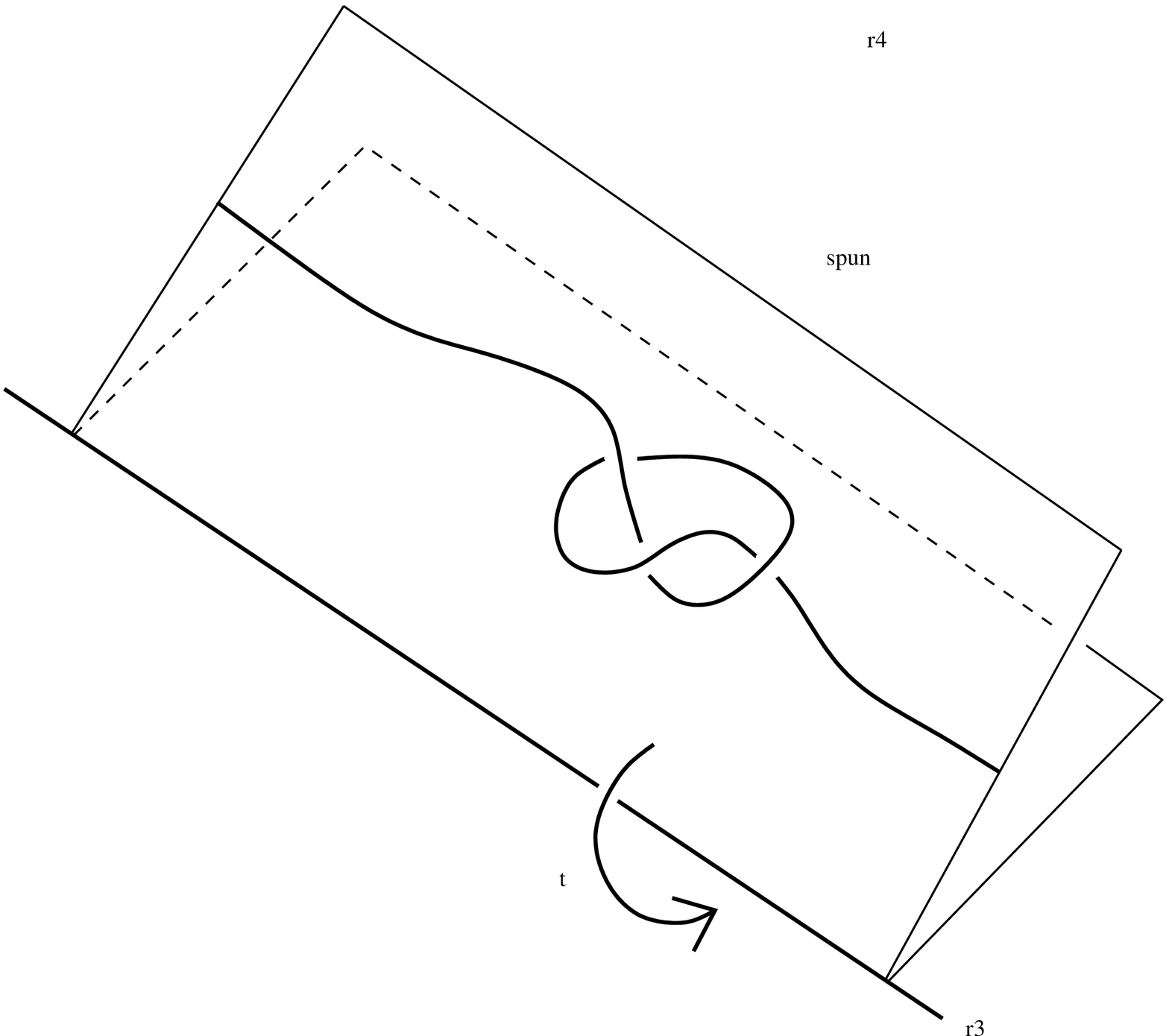}$$
\centerline{Figure 18}
}
\end{figure}
In the time between this article being accepted and published,
some progress has been made on this problem.  It turns out that,
provided $2n-3j-3 \geq 0$ the first non-trivial homotopy group
of $\Emb(\Real^j,\Real^n)$ is cyclic and in dimension $2n-3j-3$.
Moreover in these cases, a spinning construction 
$\Omega \Emb(\Real^j,\Real^n) \to \Emb(\Real^{j+1},\Real^{n+1})$
induces an epi-morphism on the first non-trivial homotopy groups
of the spaces. In particular, the spinning map
$\pi_2 \Emb(\Real,\Real^4) \to \pi_0 \Emb(S^3,\Real^6)$ is
an isomorphism -- both groups are infinite-cyclic in this
case \cite{family}.

\begin{ack}
I would like to thank Fred Cohen for teaching me about little cubes actions
 and pushing me to prove the existence of a little $2$-cubes action on
 spaces of long knots and `some kind of more general theorem,' which turned out
 to be Theorem \ref{littlecthm}.  Allen Hatcher's visit to Rochester in the
spring of 2003, and his preprint \cite{Hatcher4} on the homotopy type of $\K$,
 were of immense help when it came to formulating Theorem \ref{freeness}. 
I would like to thank Dev Sinha,
 whose work \cite{Dev} on spaces of knots has been highly inspiring. Allen Hatcher,
 Dev Sinha, and Victor Turchin's comments on the first few iterations
 of this paper were immensely valuable. I would like to thank the Mathematics
 Department at the University of Rochester for their hospitality during my visit. 
\end{ack}


\begin{thebibliography}{999}

\bibitem{Birman}
J.~Birman, 
\emph{Braids, links, and mapping class groups,} 
Annals of Mathematics Studies, No. \textbf{82}. Princeton University Press, Princeton, N.J.;
 University of Tokyo Press, Tokyo, 1974.

\bibitem{BV}
J.M.~Boardman, R.M.~Vogt. \emph{Homotopy-everything H-spaces,}
Bull. Amer. Math. Soc. {\bf 74} (1968), 1117--1122.
 
\bibitem{Borel}
A.~Borel,
\emph{Seminar on transformation groups,} 
Ann. Math. Stud. {\bf 46} Princeton University Press, Princeton, NJ. (1960)

\bibitem{Bott}
 R.~Bott, C.~Taubes, 
 \emph{On the self-linking of knots,} 
 J. Math. Phys. {\bf 35},  5247--5287 (1994)

\bibitem{Bud}
 R.~Budney, J.~Conant, K.~Scannell, D.~Sinha.
 \emph{New perspectives on self-linking,} 
 Advances in Mathematics. {\bf 191} (2005) 78--113.

\bibitem{BudCoh}
 R.~Budney, F.R.~Cohen, 
\emph{On the homology of the space of long knots in $\Real^3$}, 
preprint. 

\bibitem{bjsj}
R.~Budney,
\emph{JSJ-decompositions of knot and link complements in $S^3$},
to appear in l'Enseignement Mathematique.

\bibitem{topknot}
R.~Budney,
\emph{Topology of spaces of knots in dimension $3$,}
preprint.

\bibitem{family}
R.~Budney,
\emph{A family of embedding spaces,}
preprint.

\bibitem{BL}
D.~Burghelea, R.~Lashof, 
\emph{The homotopy type of spaces of diffeomorphisms. I,}
 Trans. Amer. Math. Soc. \textbf{196} (1974) 1--36.

\bibitem{Cat1}
A.~Cattaneo, P.~Cotta-Ramusino, R.~Longoni. 
\emph{Configuration spaces and Vassiliev classes in any dimension,}
Algebr. and Geom. Top. {\bf 2} (2002), 949--1000.

\bibitem{Cat2}
A.~Cattaneo, P.~Cotta-Ramusino, R.~Longoni. 
\emph{Algebraic structures on graph cohomology,} 
Journal of Knot Theory and Its Ramifications, Vol. 14, No. 5 (2005) 627-640.

\bibitem{EN} D.~Eisenbud, W.~Neumann, {\em Three-dimensional link
theory and invariants of plane curve singularities}, Ann. Math. Stud. {\bf 110}. (1985)

\bibitem{Good}
T.~Goodwillie, M.~Weiss, 
\emph{Embeddings from the point of view of immersion theory: Part II,}
Geom. Topol. \textbf{3} (1999), 103-118

\bibitem{Gramain}
A.~Gramain, 
\emph{Le type d'homotopie du groupe des diff\'eomorphismes d'une surface compacte,}
 Ann. Scient. \'Ec. Norm. Sup. $4^e$ serie, t. 6, 1973, p. 53 \`a 66.

\bibitem{GramainPi1}
A.~Gramain,
\emph{Sur le groupe fondamental de l'espace des noeuds,}
 Ann. Inst. Fourier, Grenoble. \textbf{27,3} (1977), 29--44. 
 
\bibitem{Pollack}
V.~Guillemin, A.~Pollack, 
\emph{Differential Topology,} Prentice-Hall, 1974.

\bibitem{Haefliger2}
A.~Haefliger,
\emph{Differentiable embeddings of $S^n$ in $S^{n+q}$ for $q>2$,}
Ann. of Math., (1966)

\bibitem{Hatcher1}
A.~Hatcher, 
\emph{Homeomorphisms of sufficiently-large $P^2$-irreducible $3$-manifolds,}
 Topology. \textbf{15} (1976)

\bibitem{Hatcher2}
A.~Hatcher, 
\emph{A proof of the Smale conjecture,} Ann. of Math. \textbf{177} (1983)

\bibitem{Hatcher3}
A.~Hatcher, 
\emph{Basic Topology of 3-Manifolds,}
\textit{[http://www.math.cornell.\-edu/\~{\hskip 0.5mm}hatcher/3M/3Mdownloads.html]}

\bibitem{Hatcher4}
A.~Hatcher, 
\emph{Topological Moduli Spaces of Knots,} 
\textit{[http://front.\-math.ucdavis.\-edu/math.GT/9909095]}

\bibitem{Hirsch}
M.W.~Hirsch, 
\emph{Differential Topology,} Springer-Verlag. (1976)

\bibitem{Ivanov1}
N.V.~Ivanov,
\emph{Diffeomorphism groups of Waldhausen manifolds,}
Research in Topology. II. Notes of LOMI scientific seminars,
V. 66 (1976), 172--176. J. Soviet Math., V. 12, No. 1 (1979), 115--118.

\bibitem{Ivanov2}
N.V.~Ivanov,
\emph{Homotopy of spaces of automorphisms of some three-dimensional manifolds,}
DAN SSSR, V. 244, No. 2 (1979), 274--277. Soviet Mathematics-Doklady,
V. 20, No. 1 (1979), 47--50.

\bibitem{JacoShalen}
W.~Jaco, P.~Shalen, 
\emph{A new decomposition theorem for $3$-manifolds,}
 Proc. Sympos. Pure Math {\bf 32} (1978) 71--84.

\bibitem{Kirby}
R.~Kirby, 
\emph{Problems in low-dimensional topology,}
Geometric topology (Athens, GA, 1993), 35--473, AMS/IP Stud. Adv. Math., 2.2, Amer.
Math. Soc., Providence, RI, 1997.

\bibitem{KS}
R.~Kirby, L.~Siebenmann, 
\emph{Foundational Essays on topological manifolds, smoothings, and triangulations,}
 Annals of math. Stud. \textbf{88} Princeton University Press. (1977)

\bibitem{Kneser}
H.~Kneser, 
\emph{Geschlossene Fl\"achen in dreidimensionalen Mannigfaltigkeiten,}
 Jahresbericht der Deut. Math. Verein., {\bf 38} (1929), 248--260.

\bibitem{komi}
H.~Kodama, P.~Michor,
\emph{The homotopy type of the space of degree 0 immersed plane curves,}
 Revista Matem\'atica Complutense {\bf 19} (2006), no. 1, 227-234.

\bibitem{Kohno}
T.~Kohno, 
\emph{Linear representations of braid groups and classical
 Yang-Baxter equations,} Cont. Math. {\bf 78} (1988), 339--363.

\bibitem{Kont}
M.~Kontsevich, 
\emph{Vassiliev's knot invariants,} I.~M.~Gelfand Seminar
 (S.~Gelfand and S.~Gindikin, eds.) Adv. Soviet. Math., {\bf 16}, Amer. Math.
 Soc., Providence, 1993, pp. 137--150.

\bibitem{Lith}
R.A.~Litherland,
\emph{Deforming twist-spun knots,}
Trans. Amer. Math. Soc. {\bf 250} (1979), 311--331.
 
\bibitem{MSS}
M.~Markl, S.~Shnider, J.~Stasheff,
\emph{Operads in algebra, topology and physics,}
AMS Mathematical surveys and monographs. Vol {\bf 96}.(2002)

\bibitem{May1}
J.P.~May, 
\emph{The Geometry of Iterated Loop Spaces,}
 Lecture Notes in Mathematics. \textbf{271} (1972)

\bibitem{May2}
 J.P.~May, 
 \emph{$E_{\infty}$ spaces, group completions, and permutative categories,}
 London Math. Soc.  Lecture Notes Series \textbf{11}, 1974, 61--93.

\bibitem{McClure}
J.~McClure, J.~Smith. 
\emph{Cosimplicial objects and little $n$-cubes. I,} 
\textit{[http://\-front.\-math.ucdavis.\-edu/math.QA/0211368]}

\bibitem{Mor}
C.~Morlet, 
\emph{Plongement et automorphismes de vari\'et\'es,} Notes multigraphi\'ees,
Coll\`ege de France, Cours Peccot (1969).

\bibitem{Neumann}
W.~Neumann,
\emph{Notes on geometry and $3$-manifolds,}
Low-dimensional Topology. Bolyai Society Mathematical Studies {\bf 8} (1999), 191--267.
 
\bibitem{Nielsen}
J.~Nielsen,
\emph{Untersuchungen zur theorie der geschlossenen zweiseitigen Fl\"achen I,}
Acta Math. {\bf 50} (1927), 189--358.

\bibitem{Pal}
R.~Palais, 
\emph{Local triviality of the restriction map for embeddings,}
Comment. Math. Helv. {\bf 34} (1960) 305--312.

\bibitem{salvatore}
P.~Salvatore,
\emph{Knots, operads and double loop spaces,}
preprint.

\bibitem{DevKevin}
K.~Scannell, D.~Sinha, 
\emph{A one-dimensional embedding complex,} 
Journal of Pure and Applied Algebra {\bf 170} (2002), No. 1, 93--107.

\bibitem{Sch}
H.~Schubert, 
\emph{Die eindeutige Zerlegbarkeit eines Knoten in Primknoten,}
 Sitzungsber. Akad. Wiss. Heidelberg, math.-nat. KI., \textbf{3:57-167}. (1949)

\bibitem{Dev}
D.~Sinha, 
\emph{The topology of spaces of knots,} preprint,
\textit{[http://\-front.math.ucdavis.\-edu/math.AT/0202287]}

\bibitem{Dev2}
D.Sinha, 
\emph{Operads and knot spaces.} 
J. Amer. Math. Soc. {\bf 19} (2006), no. 2, 461--486

\bibitem{Smale}
S.~Smale, 
\emph{ Diffeomorphisms of the $2$-sphere,}
Proc. Amer. Math. Soc. \textbf{10} (1959) 621--626.

\bibitem{Steenrod}
N.~Steenrod, 
\emph{A convenient category of topological spaces,}
 Mich. Math. J. \textbf{14} no. 2, (1967) 133--152.

\bibitem{Thurston}
W.~Thurston, 
\emph{Hyperbolic structures on 3-manifolds, I. Deformation of acylindrical manifolds,}
 Ann. of Math. \textbf{124} (1986), 203--246.

\bibitem{Tour}
V.~Turchin, 
\emph{On the homology of the spaces of long knots,} 
Advances in Topological Quantum Field Theory. 
NATO Sciences series by Kluwer 2005, pp 23-52.

\bibitem{Vass}
V.~Vassiliev, 
\emph{Complements of discriminants of smooth maps: topology  and applications,}
  Translations of Mathematical Monographs, {\bf 98}.
 American Mathematical Society, Providence, RI, 1992.

\bibitem{Volic}
 I.~Volic, 
 \emph{Finite-type knot invariants and calculus of functors,}
  Brown University. (2003)

\bibitem{Zie}
H.~Zieschang, E.~Vogt, H.~Coldeway,
\emph{Surfaces and planar discontinuous groups,}
Lecture Notes in Math., no. 835, Springer, 1980.



\end{thebibliography}
\end{document}